\documentclass[11pt]{article} 
\usepackage[utf8]{inputenc} 

\usepackage[margin=1.25in]{geometry} 
\usepackage{graphicx} 

\usepackage{amssymb}
\usepackage{amsthm}
\usepackage{amsmath,amsfonts,amssymb}

\usepackage[usenames,dvipsnames]{xcolor}
\usepackage[colorlinks=true, pdfstartview=FitV, linkcolor=blue, citecolor=blue, urlcolor=blue]{hyperref}

\newcommand{\de}{\delta}
\newcommand{\la}{\lambda}
\newcommand{\CC}{\mathbb{C}}
\newcommand{\RR}{\mathbb{R}}
\newcommand{\ringphi}{\mathring{\phi}}
\newcommand{\ringpsi}{\mathring{\psi}}
\newcommand{\sh}{\mathcal{H}}

\numberwithin{equation}{section}

\newtheorem{theorem}{Theorem}[section]
\newtheorem{remark}[theorem]{Remark}
\newtheorem{proposition}[theorem]{Proposition}

\newcommand{\ringf}{\mathring{f}}
\newcommand{\ringg}{\mathring{g}}
\newcommand{\eps}{\varepsilon}

\newcommand{\e}{\varepsilon}
\newcommand{\ENZ}{\Omega \setminus \overline{D}}
\newcommand*{\R}{\ensuremath{\mathbb{R}}}
\newcommand{\epsenz}{\varepsilon_{\rm \scriptscriptstyle ENZ}}

\let\originalleft\left
\let\originalright\right
\renewcommand{\left}{\mathopen{}\mathclose\bgroup\originalleft}
\renewcommand{\right}{\aftergroup\egroup\originalright}

\renewcommand*{\tilde}{\widetilde}
\renewcommand*{\hat}{\widehat}

\usepackage{titlesec}

\newcommand{\addperiod}[1]{#1.}
\titleformat{\section}
   {\normalfont\Large}{\thesection.}{0.5em}{}
\titleformat*{\subsection}{\normalfont\large}
\titleformat{\subsubsection}[runin]
  {\bfseries}
  {\thesubsubsection.}
  {0.5em}
  {\addperiod}
\titleformat*{\subsubsection}{\bfseries}
\titleformat*{\paragraph}{\bfseries}
\titleformat*{\subparagraph}{\large\bfseries}

\title{\bf \Large Transverse Magnetic ENZ Resonators: Robustness and Optimal Shape Design}

\author{Robert V. Kohn\textsuperscript{1}
\thanks{Courant Institute of Mathematical Sciences, New York University. Email:
{\footnotesize \href{mailto:kohn@cims.nyu.edu}{kohn@cims.nyu.edu}.}
}
\and
Raghavendra Venkatraman\textsuperscript{1}
\thanks{Courant Institute of Mathematical Sciences, New York University. Email:
{\footnotesize \href{mailto:raghav@cims.nyu.edu}{raghav@cims.nyu.edu}.}
}
}
\date{\today}

\begin{document}

\maketitle

\begin{abstract}
We study certain ``geometric-invariant resonant cavities'' introduced by Liberal,
Mahmoud, and Engheta in a 2016 Nature Communications paper. They are cylindrical
devices modeled using the transverse magnetic reduction of Maxwell's equations,
so the mathematics is two-dimensional. The cross-section consists of a
dielectric inclusion surrounded by an ``epsilon-near-zero'' (ENZ) shell. When
the shell has just the right area, its interaction with the
inclusion produces a resonance. Mathematically, the resonance is
a nontrivial solution of a 2D divergence-form Helmoltz equation
$\nabla \cdot \left( \eps^{-1}(x,\omega) \nabla u \right) + \omega^2 \mu u = 0$,
where $\eps(x,\omega)$ is the (complex-valued) dielectric permittivity, $\omega$ is the frequency,
$\mu$ is the magnetic permeability, and a homogeneous Neumann condition is imposed at the
outer boundary of the shell. This is a nonlinear eigenvalue problem, since $\eps$ depends on $\omega$.
Use of an ENZ material in the shell
means that $\eps(x,\omega)$ is nearly zero there, so the PDE is rather singular. Working
with a Lorentz model for the dispersion of the ENZ material, we put the discussion of
Liberal et.~al.~on a sound foundation by proving the existence of the anticipated resonance when
the loss parameter of the Lorentz model is sufficiently small. Our analysis is perturbative
in character, using the implicit function theorem despite the apparently singular form of
the PDE. While the existence of the resonance depends only on the area of the ENZ
shell, its quality (that is, the rate at which the resonance decays) depends on the
shape of the shell. It is therefore natural to
consider an associated optimal design problem: what shape shell gives the
slowest-decaying resonance? We prove that if the dielectric inclusion is
a ball then the optimal shell is a concentric annulus. For an inclusion of any
shape, we study a convex relaxation of the design problem using tools
from convex duality. Finally, we discuss the conjecture that our relaxed problem amounts
to considering homogenization-like limits of nearly optimal designs.
\end{abstract}
\footnotetext[1]{The authors warmly thank Nader Engheta for bringing this problem to our attention, and for many helpful
discussions. Both authors gratefully acknowledge support from the Simons Foundation through its Collaboration
on Extreme Wave Phenomena (grant 733694). RVK also gratefully acknowledges support from the National Science
Foundation (grant DMS-2009746).}
\newpage

\setcounter{tocdepth}{2}
\tableofcontents

\section {Introduction} \label{sec:introduction}

This paper is motivated by a 2016 article by Liberal et al, which discusses how
epsilon-near-zero (ENZ) materials can be used to design ``geometry-invariant resonant
cavities'' \cite{liberal2016geometry} .
We focus on a class of examples involving the
transverse magnetic reduction of the time-harmonic Maxwell system, obtained by taking
$H = (0,0,u(x_1,x_2))$ and $E =\frac{1}{i \omega \eps}(-\partial_2 u, \partial_1 u, 0)$
in
\begin{equation} \label{time-harmonic maxwell}
	\nabla \times H  = - i \omega \eps(x) E, \quad \nabla \times E = i \omega \mu(x) H .
\end{equation}
Thus we shall be working with the Helmholtz equation
\begin{equation} \label{helmholtz}
	\nabla \cdot \left( \frac{1}{\eps(x)} \nabla u \right) + \omega^2 \mu(x) u = 0
\end{equation}
in a bounded two-dimensional domain $\Omega$. Here $\omega$ is the frequency, and
$\eps = \eps(x_1, x_2)$, $\mu = \mu(x_1,x_2)$ are the dielectric permittivity and
magnetic permeability at this frequency.

It is easy to see from \eqref{helmholtz} what is special about ENZ materials in
the transverse-magnetic setting. Indeed, if $\eps(x)$ is near zero in some ``ENZ region,''
then $\frac{1}{\eps(x)} \nabla u$ can avoid being large only by $\nabla u$ being small in this
region. In the limit as $\eps \rightarrow 0$ in the ENZ region, we are not solving a PDE there
but rather choosing a constant value for $u$. While the
solution of a PDE depends sensitively on its domain and coefficients, the constant value of $u$ in the
ENZ region should be much less sensitive. In fact, in many settings it is only the \emph{area} of the ENZ region
that matters (to leading order, as $\eps \rightarrow 0$ in the ENZ region). This effect has been used, for
example, to design entirely new types of waveguides \cite{silveirinha2006tunneling,silveirinha2007design,silveirinha2007theory}; for recent reviews of these and other applications, see \cite{liberal2017rise,niu2018epsilon}.

\begin{figure}[h]
	\centerline{\includegraphics[scale=.7]{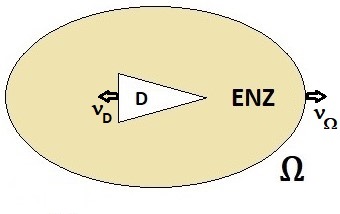}}
	\caption{Our resonator consists of a region $D$ in which $\eps = \eps_0 \eps_D(\omega)$
		surrounded by a shell $\ENZ$ in which $\eps = \eps_0 \epsenz (\omega)$. The material
		in the shell is ENZ, whereas $D$ is occupied by an ordinary dielectric.}
	\label{fig:resonator-geometry}
\end{figure}

We now specify more precisely the PDE problem considered in this paper. The Helmholtz equation
\eqref{helmholtz} will be solved in a bounded domain $\Omega \subset \R^2$ with a core-and-shell
structure: it consists of a region $D \subset \Omega$ containing an ordinary dielectric, surrounded by a
shell $\Omega \setminus D$ containing an ENZ material (see Figure \ref{fig:resonator-geometry}). At the outer
boundary $\partial \Omega$ we take $\partial u / \partial \nu_\Omega = 0$. (In the language of the underlying
Maxwell system, the outer boundary is a perfect electric conductor.) Following \cite{liberal2016geometry},
we shall ignore the spatial and frequency dependence of $\mu$, since it is negligible in the intended applications;
thus we set $\mu(x) = \mu_0$ to be the permeability of free space. The dielectric permittivity $\eps$ is
constant in each material:
\begin{equation} \label{eps-vs-relative-permeabilities}
	\eps(x) = \left\{ \begin{array}{ll}
		\eps_0 \, \eps_D (\omega) & \mbox{in $D$}\\
		\eps_0 \, \epsenz (\omega) & \mbox{in $\ENZ$}
	\end{array}
	\right.
\end{equation}
where $\eps_0$ is the permittivity of free space. While our method is more general, we shall use
a \emph{Lorentz model} for the ENZ material:
\begin{equation} \label{lorentz-model}
\epsenz (\omega) = \eps_\infty \left( 1 + \frac{\omega_p^2}{\omega_0^2 - \omega^2 - i \omega \gamma} \right)
\end{equation}
Here $\eps_\infty$, $\omega_0$, $\omega_p$ and $\gamma$ are nonnegative real numbers.
Notice that in the lossless case $\gamma = 0$
this function vanishes when $\omega = \omega_* := \sqrt{\omega_p^2 + \omega_0^2}$; the resonant frequency of our ENZ-based
resonator will be very near this ``ENZ frequency.'' As for $\eps_D$:
following \cite{liberal2016geometry}, we shall ignore losses there, taking
$\eps_D (\omega)$ to be real-valued and positive for real-valued $\omega$ near $\omega_p$.

With the preceding conventions, and writing $c = 1/\sqrt{\eps_0 \mu_0}$ for the speed of light,
the Helmholtz equation \eqref{helmholtz} becomes
\begin{eqnarray}
	\nabla \cdot \left( \frac{1}{\eps_D (\omega)} \nabla u \right) + \omega^2 c^{-2} u = 0 && \quad \mbox{in $D$ and}
	\label{helmholtz-D}\\
	\nabla \cdot \left( \frac{1}{\epsenz (\omega)} \nabla u \right) + \omega^2 c^{-2} u = 0 && \quad \mbox{in $\ENZ$}
	\label{helmholtz-ENZ}
\end{eqnarray}
with the understanding that $u$ and $\frac{1}{\eps} \nabla u \cdot \nu_D$ are continuous at $\partial D$, and
that $u$ satisfies $\partial u / \partial \nu_\Omega = 0$ at $\partial \Omega$. Following \cite{liberal2016geometry} ,
we shall refer to
a nontrivial solution as a \emph{resonance}. It should be noted, however, that when the loss parameter $\gamma$ is
positive $\epsenz$ is complex, so the resonance $u$ and the resonant frequency $\omega$ will also be
complex. Since the time-harmonic Maxwell equations are obtained by considering electric and magnetic fields
of the form $e^{-i \omega t} E(x)$ and $e^{-i \omega t} H(x)$, and since physical solutions should decay in time,
we expect (and we will find) that in the presence of loss, the imaginary part of $\omega$ is negative.

While the preceding discussion is accurate, it ignores an important feature of our analysis. Indeed, in
discussing the Lorentz model \eqref{lorentz-model} we wrote $\epsenz = \epsenz (\omega)$, treating
the constants $\eps_\infty$, $\omega_0$, $\omega_p$, and $\gamma$ as being fixed. Actually, the dependence of
$\epsenz$ on the loss parameter $\gamma$ is crucial to our analysis. In fact, when we prove the
existence of a resonance in Section \ref{sec:accounting-for-dispersion-and-loss}, our main tool is
perturbation theory in $\gamma$, and the resonant frequency is a function of $\gamma$. When the dependence
of $\epsenz$ on $\gamma$ is important, we shall write $\epsenz(\omega,\gamma)$ rather than $\epsenz(\omega)$
for the Lorentz model \eqref{lorentz-model} (see for example equation \eqref{suffices-to-find} below).

Equations \eqref{helmholtz-D}--\eqref{helmholtz-ENZ} are not a conventional eigenvalue problem, since $\eps_D$ and
$\epsenz$ depend on $\omega$. The fundamental insight of \cite{liberal2016geometry} in this context was that \emph{there should
	nevertheless be a solution near the ENZ frequency when $\gamma$ is small, provided that the area of the ENZ region
	satisfies a certain consistency condition}. This result is interesting (and potentially useful) because the consistency condition
involves \emph{only} the area of the ENZ region. Usually, for a PDE, the existence of a resonance at a given frequency
depends sensitively on the geometry of the domain. Our ENZ-based resonator is different: it has a resonance near the
ENZ frequency \emph{regardless} of the shape of the ENZ region, provided only that the area of this region is right. Thus,
the design of resonators with a given resonant frequency becomes easy (given the availability of a material with $\eps \approx 0$
at that frequency).

The paper \cite{liberal2016geometry}  uses physical insight to find the condition on the area of the ENZ region, and it uses
numerical simulations to confirm that the anticipated resonances exist in many examples. Our work complements its
contributions by proving the existence of such resonances and studying their dependence on the ENZ material's loss parameter
$\gamma$. In particular, we provide a rather complete understanding about how the geometry of the ENZ region influences the
rate at which the resonance decays. It is natural to ask how the shape of the ENZ shell can be chosen to minimize the decay rate.
When $D$ is a disk, we show that the optimal ENZ shell is a concentric annulus; for more general $D$, a similarly
explicit solution is probably not possible, but we are nevertheless able to estimate the optimal decay rate by
considering a certain convex optimization.

Our account has thus far emphasized the physical character of the problem. To communicate the mathematical character
of our work, it convenient (and indeed necessary) to consider what happens when we \emph{ignore} the frequency-dependence
of $\eps_D$ and $\epsenz$. After multiplying both equations \eqref{helmholtz-D}--\eqref{helmholtz-ENZ} by $\eps_D$,
our PDE \eqref{helmholtz} takes the form
\begin{align} \label{simplified-pde}
	\nabla \cdot \e_\delta^{-1}\nabla  u + \lambda u &= 0 \quad \mbox{ in } \Omega\\
	\partial_{\nu_\Omega} u &= 0 \quad \mbox{ at } \partial \Omega\, , \nonumber
\end{align}
with the conventions that $\lambda = \omega^2 c^{-2} \eps_D$ and $\delta = \epsenz /\eps_D$, and the notation
\begin{equation} \label{simplified-eps}
	\eps_\delta(x) := \left\{
	\begin{array}{cc}
		\delta & x\in \ENZ \\
		1      & x \in D\, .
	\end{array}
	\right.
\end{equation}
Since $\delta$ is just a parameter, this \emph{is} a linear eigenvalue problem. It appears to be rather
singular in its dependence on $\delta$, since the PDE in $\ENZ$ is now
$\nabla \cdot (\delta^{-1} \nabla u) + \lambda u = 0$ and we are interested $\delta$ near $0$. But
it can be desingularized by a suitable ansatz, as we shall explain in Section \ref{sec:analysis-without-dispersion}.
(In $\ENZ$ the ansatz has $u = 1 + \delta f (x)$, so that $ \delta^{-1} \nabla u = \nabla f$ is no longer singular.)

When $\delta$ is real and positive, it is a basic result about second-order elliptic PDE that \eqref{simplified-pde}
can have a nontrivial solution (a resonance) only for a discrete set of $\lambda$'s, which must be nonnegative. Our
main result about \eqref{simplified-pde}, Theorem \ref{t.main1}, identifies the (infinite) set of $\lambda$'s for which
such a result holds even for \emph{complex-valued} $\delta$ in a neighborhood of $0$; moreover it shows that for
each such resonance, $u = u_\delta$ and $\lambda = \lambda_\delta$ are complex analytic in their dependence on $\delta$.
(Our results agree with those in \cite{liberal2016geometry}  concerning the possible values
of $\lambda_0 = \lim_{\delta \rightarrow 0} \lambda_\delta$ and $u_0(x) = \lim_{\delta \rightarrow 0} u_\delta (x)$.)
Besides proving analyticity, our work gives easy access to the Taylor expansions of $u_\delta$ and $\lambda_\delta$; in
particular, it identifies the asymptotic electric field in the ENZ region (in other words, the limiting value of
$\delta^{-1} \nabla u$ in $\ENZ$ as $\delta \rightarrow 0$), and it shows how the shape of the ENZ shell determines
the leading-order correction to $\lambda_0$ when $\delta \neq 0$ (that is, the value of
$\lambda'(0)$).

Let us say a word about the proof of Theorem \ref{t.main1}. The arguments draw inspiration from those used to study perturbations
of eigenvalues in more standard settings. Due to the singular character of our operator, however, we must solve
PDE's in $D$ and $\ENZ$ \emph{separately}, rather than ever solving \eqref{simplified-pde} in the entire domain $\Omega$.
Our analysis begins by showing how the Taylor expansions of $u_\delta$ and $\lambda_\delta$ can be determined
term-by-term. While analyticity (with respect to $\delta$) can be proved by majorizing the resulting expansion, we pursue
a different approach -- demonstrating analyticity by an application of the implicit function theorem.

We are not the first to consider operators of the form
$\nabla \cdot ( \eps_\delta^{-1} \nabla u)$ in which $\eps_\delta$
takes only the values $1$ and $\delta$ and the focus is on behavior near
$\delta = 0$. This operator and others like it arise, in particular, when considering the effective
behavior or band structure of high-contrast composites
\cite{Bru,Chen-Lipton-ARMA2013,Chen-Lipton-MMS2013,Hempel-Lienau-2000,Fortes-Lipton-Shipman-2010,Fortes-Lipton-Shipman-2011}.
Our treatment of \eqref{simplified-pde} has some features in common with the work just cited,
as we discuss in more detail near the end of this section.

Returning to the physical problem with dispersion and loss, \eqref{helmholtz-D}--\eqref{helmholtz-ENZ}: the existence of
resonances and their analytic dependence on the loss parameter $\gamma$ follows easily from Theorem \ref{t.main1} by an
application of the implicit function theorem. Indeed, it suffices to find a complex-valued function $\omega(\gamma)$
such that
\begin{equation} \label{suffices-to-find}
\lambda_{\epsenz (\omega(\gamma),\gamma)/\eps_D(\omega(\gamma))} = \omega^2(\gamma) c^{-2} \eps_D(\omega(\gamma))
\end{equation}
and such that $\omega(0)$ is the ENZ frequency (the one where $\epsenz$ vanishes when $\gamma = 0$).
We show in Theorem \ref{t.result-with-dispersion} that the implicit function theorem is applicable, and that the leading-order
dependence of $\omega(\gamma)$ (that is, $\omega'(0)$) depends on the geometry of the ENZ region only
through $\lambda'(0)$. Also of note: we show that $\omega'(0)$ is purely imaginary. It follows that
the frequency where the resonance occurs (the real part of $\omega(\gamma)$) differs very little from
the ENZ frequency (the difference is at most of order $\gamma^2$).

For an experimentalist creating a resonator using the framework of this paper,
the ENZ material to be used in $\ENZ$ would typically be fixed, and
therefore the ENZ frequency $\omega_*$ (defined by \eqref{enz-frequency})
would also be fixed. As mentioned earlier, the experimentalist's choices of $D$ and $\Omega$
must satisfy a consistency condition that depends on $\omega_*$. This is
discussed in Section \ref{subsubsec:consistency}, but we summarize the main
impact here: (i) while the shape of $D$ is unconstrained, its size must satisfy a
certain (open) condition; (ii) once both the ENZ frequency and $D$ are fixed,
the consistency condition constrains $\Omega$ only by fixing the area of the
ENZ region $\ENZ$.

It is natural to ask how the shape of the ENZ region should be chosen to optimize the associated resonance. Since
the imaginary part of $\omega$ gives the rate at which the resonance decays, this amounts to asking what shape
minimizes $|\omega'(0)|$. The analysis just summarized reduces this to asking what shape minimizes
$|\lambda'(0)|$. Our results on the function $\lambda(\delta)$ include a variational characterization of this number: it is
a constant times
\begin{equation} \label{early-var-prin}
\min_\phi \int_{\ENZ} \, \frac{1}{2} | \nabla \phi|^2 - \lambda_0 \phi \, dx + \int_{\partial D}  f  \phi \, d\sh^1
\end{equation}
for a particular choice of the function $f$ (see \eqref{variational-characterization}). Since the minimum value of
\eqref{early-var-prin} is negative, our optimal design problem takes the form
\begin{equation} \label{optimal-design-pbm-ver1}
\max_\Omega \min_\phi \int_{\ENZ} \, \frac{1}{2} | \nabla \phi|^2 - \lambda_0 \phi \, dx + \int_{\partial D}  f  \phi \, d\sh^1 \, ,
\end{equation}
with the understanding that $\Omega$ varies over domains that contain $D$ and remain within some fixed region $B$. The
domain $\Omega$ can be represented by a function $\chi(x)$, defined on $B \setminus \overline{D}$, which takes the
value $1$ on $\ENZ$ and $0$ outside $\Omega$. With this convention, \eqref{optimal-design-pbm-ver1} becomes
$$
\max_{\chi(x) \in \{0,1\} }  \min_\phi \int_{B \setminus \overline{D}}
\chi(x) \left( \frac{1}{2} |\nabla \phi|^2 - \lambda_0 \phi \right) \, dx + \int_{\partial D}  f  \phi \, d\sh^1 .
$$
In the language of structural optimization (see e.g. \cite{Allaire}), this is a
\emph{compliance optimization} problem. Such problems are well-understood for mixtures of two
nondegenerate materials (that is, when $\chi(x)$ takes two values that are both positive).  In the present
more degenerate setting, methods from homogenization cannot be applied directly. But taking inspiration
from that theory, we show in Section \ref{sec:the-optimal-design-problem} that the value of
\eqref{optimal-design-pbm-ver1} is upper-bounded by the value of the simple-looking convex optimization
\begin{equation} \label{optimal-design-pbm-dual}
\min_\phi \int_{B \setminus \overline{D}} \left( \frac{1}{2} |\nabla \phi |^2 - \lambda_0 \phi \right)_+ \, dx +
\int_{\partial D}  f  \phi \, d\sh^1 \, .
\end{equation}
Moreover, we argue (though we do not prove) that this bound is actually sharp. When the inclusion $D$ is a ball
we can say much more: the optimal $\Omega$ is in fact a concentric ball (and the upper
bound \eqref{optimal-design-pbm-dual} is indeed sharp in this case).

Let us briefly discuss some related work.

\begin{itemize}
	
\item The physics literature includes many papers on devices made using ENZ materials, including
quite a few that can be modeled by Helmholtz equations like \eqref{helmholtz}. Many of these papers
raise issues comparable to those considered in the present work. Our recent paper \cite{KV1}
studied a phenomenon known as photonic doping; it provided a mathematical foundation for and an
improved understanding of an application of ENZ materials considered in \cite{liberal-etal-photonic-doping}
(see also \cite{silveirinha2007design}). That work
involved \emph{scattering}, whereas the present work involves \emph{resonance}. Therefore the analysis
in this paper is substantially different from that of \cite{KV1}, though there are of course some parallels.
In particular, Section \ref{sec:analysis-without-dispersion} of this paper shows how the perturbation theory
of eigenvalue problems can be adapted to the ENZ setting, while our earlier paper was concerned instead
with the perturbation theory of boundary-value problems.

\item As already mentioned earlier, our treatment of \eqref{simplified-pde} has some features in common with
\cite{Bru,Chen-Lipton-ARMA2013,Chen-Lipton-MMS2013,Hempel-Lienau-2000,Fortes-Lipton-Shipman-2010,Fortes-Lipton-Shipman-2011,KV1}.
Preparing to say more on this topic, we remind the reader that
\eqref{simplified-pde} is an eigenvalue problem for the
divergence-form operator $\nabla \cdot \left( \e_\delta^{-1}\nabla  u \right)$,
whose coefficient $\e_\delta$ is piecewise constant (equal to $1$ in $D$
and $\delta$ in $\ENZ$). We show in Section \ref{sec:analysis-without-dispersion}
that the quantities of interest are complex analytic functions of $\delta$
near $\delta=0$.
\medskip

For any function of a complex variable $\delta$, there are two rather distinct
approaches to proving its analyticity. One is to show that the function
is complex differentiable in $\delta$; this is what we do in
Section \ref{sec:analysis-without-dispersion}. The other is to identify
the function's Taylor expansion then prove its convergence; this is the
approach taken in
\cite{Bru,Chen-Lipton-ARMA2013,Chen-Lipton-MMS2013,Fortes-Lipton-Shipman-2010,Fortes-Lipton-Shipman-2011,KV1},
which consider problems closely related to
ours. Among these references, the papers by Chen \& Lipton and Fortes, Lipton \& Shipman
have the strongest connections to our setting, since they too consider spectral problems.
This work studies the band structure of certain periodic high-contrast composites; thus
its physical motivation is quite different from ours. However, the PDEs
considered in these papers are closely analogous to \eqref{simplified-pde}
(except for being solved on a period cell, with a Bloch boundary condition).
Therefore it is not surprising that in these papers, as in Section \ref{sec:analysis-without-dispersion}, one finds each
successive term of the Taylor expansion by considering (separate) PDE problems in two complementary
material regions; moreover, the Taylor expansions found in these papers have a character quite similar
to ours. (Since the work just discussed concerns the band structure of periodic high-contrast
composites, let us also mention an earlier paper \cite{Hempel-Lienau-2000}, which
achieves impressive insight by means other than Taylor expansion.)
\medskip

Rather than majorize the Taylor expansion, our proof of
analyticity uses the implicit function theorem to show that the eigenvalue $\lambda_\delta$
and the (suitably normalized) eigenfunction $u_\delta$ of \eqref{simplified-pde} are complex-differentiable
functions of $\delta$ near $\delta = 0$. The fact that perturbation theory for (simple) eigenvalues can
be done using the implicit function theorem has been understood at least since 1955
\cite{Rosenbloom-1955}. While this approach does not immediately give a radius of analyticity,
extensions of that type have been discussed in some settings \cite{Kloeckner-2019}.

\item As we explain in Section \ref{subsec:relaxation}, the passage from
\eqref{optimal-design-pbm-ver1} to \eqref{optimal-design-pbm-dual}
involves considering the possibility that the optimal $\Omega$ is a
homogenization limit of domains with many small holes. Our optimal
design problem can be regularized by including a penalty term involving the
\emph{perimeter} of $\ENZ$. It is known that inclusion of such a penalty prevents
homogenization (see e.g.~\cite{bourdin-chambolle}). However, if the unpenalized
optimization requires homogenization then the solution of the penalized problem will depend
strongly on the presence and strength of the penalization. Therefore we do
not consider the use of perimeter penalization in the present work.
\end{itemize}

We close this Introduction by summarizing the organization of the paper. Section \ref{sec:analysis-without-dispersion} contains our
study of the PDE \eqref{simplified-pde}. It is the longest section, since much of our success lies in finding a convenient way to
desingularize the problem. Section \ref{sec:accounting-for-dispersion-and-loss} combines the results of Section \ref{sec:analysis-without-dispersion}
with the implicit function theorem to show the existence of a resonance near the ENZ frequency, and to
consider the effects of dispersion and loss. Finally, Section \ref{sec:the-optimal-design-problem} presents our results on the optimal design
problem (choosing the shape of the ENZ region to minimize the effect of loss).

\section{Analysis without dispersion}
\label{sec:analysis-without-dispersion}

In this section we study the eigenvalue problem \eqref{simplified-pde}. Our main result is
the existence of an eigenfunction $u_\de$ with eigenvalue $\lambda_\de$, both depending
complex-analytically on $\de$ in a neighborhood of $0$, provided that
$\lambda_0 = \lim_{\de \rightarrow 0} \lambda_\de$ satisfies a certain consistency
condition. Our proof shows in addition that
$\lambda_\de$ is a \emph{simple} eigenvalue, in other words its eigenspace is one-dimensional.

We start, in Section \ref{subsec:preliminaries-etc}, with some preliminaries and a full statement of the result.
Then we show, in Sections \ref{subsec:leading-order-terms} -- \ref{subsec:higher-order-terms},
how the Taylor expansions of $u_\de$ and $\lambda_\de$ can be determined term-by-term if one
assumes analyticity. Finally, in Section \ref{subsec:analyticity} we use the implicit function theorem to
prove the existence of $u_\de$ and $\lambda_\de$ depending analytically on $\de$.

Our analysis shows, roughly speaking, that the perturbation theory of eigenvalues and eigenfunctions
can be applied to the singular-looking operator $\nabla \cdot \e_\de^{-1}(x)\nabla$ with estimates that
are uniform in $\de$ (and that $\delta = 0$ is a removable singularity).

\subsection{Preliminaries and a statement of our result about $u_\de$ and $\lambda_\de$}
\label{subsec:preliminaries-etc}

As discussed in the Introduction,
we are interested throughout this paper in a bounded two-dimensional domain
$\Omega$ with a subset $D$ (see Figure \ref{fig:resonator-geometry}).
Both domains are assumed to be Lipschitz (that is, their boundaries are locally the graphs of Lipschitz
functions) and connected, and $\overline{D}$ does not touch $\partial \Omega$. We also assume that
$D$ is simply connected, so that the ``ENZ region'' $\ENZ$ is a connected set which can be viewed as a
shell surrounding $D$. While the $\Omega$ shown in Figure \ref{fig:resonator-geometry} is simply
connected, we do not assume this; rather, $\ENZ$ can have one or more holes -- in which case the
boundary condition $\partial_{\nu_\Omega} u = 0$ in \eqref{simplified-pde} applies at the boundary
of each hole.

\subsubsection{The function $\psi_{d,\lambda_0}$; a normalization; and the consistency condition}
\label{subsubsec:psi-d}

It is natural to begin by finding the Taylor expansion of $u_\delta$ and $\lambda_\delta$,
assuming existence and analyticity. As a reminder, our goal is to solve
\begin{equation} \label{e.evp}
\begin{aligned}
\nabla \cdot \tfrac{1}{\e_\delta} \nabla u_\de + \lambda_\de u_\de &= 0 \quad \mbox{in $\Omega$}\\
\partial_{\nu_\Omega} u_\de &= 0 \quad \mbox{at $\partial \Omega$} \, ,
\end{aligned}
\end{equation}
where $\e_\delta (x) := 1$ for $x \in D$ and $\e_\delta(x) = \delta$ for $x \in \ENZ$. Proceeding formally for the
moment, we seek a solution of the form
\begin{align} \label{e.lambdaexp}
	\lambda_\delta = \lambda_0 + \delta \lambda_1 + \delta^2 \lambda_2 + \ldots
\end{align}
and
\begin{equation} \label{e.ansatz}
\begin{aligned}
u_\delta(x) := \left\{ \begin{array}{ll}
1 + \delta \phi_1 + \delta^2 \phi_2 + \ldots  & \mbox{ if } x \in \Omega \setminus \overline{D} \\
\psi_0 + \delta \psi_1 + \delta^2 \psi_2 + \ldots & \mbox{ if } x \in D \, .
\end{array}
\right.
\end{aligned}
\end{equation}
This should, of course, not be possible for all choices of $\lambda_0$; the condition that determines the permissible values of
$\lambda_0$ will be given presently (see \eqref{consistency-first-version}).

The expansion of $u_\delta$ begins with $1$ in $\ENZ$ because (as discussed in the Introduction) we expect $u_\delta$ to be
constant to leading order in $\ENZ$. There is no loss of generality taking the leading-order constant to be $1$, since
multiplying an eigenfunction by a constant gives another eigenfunction. But this normalization only affects the leading-order
term, whereas to fix $u_\delta$ we need a condition that applies to all orders in $\delta$. It is convenient to use the
normalization
\begin{equation} \label{normalization-using-psi0}
\int_\Omega u_\de (x)u_0(x) \, dx = \int_\Omega u_0^2\,dx = |\ENZ| + \int_D \psi_0^2\,dx \, ,
\end{equation}
where $u_0 = \lim_{\delta \rightarrow 0} u_\delta$ denotes the leading-order term of \eqref{e.ansatz},
\begin{equation} \label{u-zero}
u_0(x) = \left\{
	\begin{array}{cc}
		1 &  x \in \ENZ \\
		\psi_0 & x \in D \, .
	\end{array}
	\right.
\end{equation}

When we substitute the expansions \eqref{e.lambdaexp}--\eqref{e.ansatz} into the PDE \eqref{e.evp} and
focus on the leading-order behavior in $D$, we see that $\psi_0$ must solve a Helmholtz equation in $D$ with
the Dirichlet boundary condition $\psi_0 = 1$ at $\partial D$. The solution of this boundary value problem
also played a central role in our recent study of photonic doping \cite{KV1}. To emphasize the connections
between that work and this one we will use similar notation here, calling its solution $\psi_{d, \lambda_0}$.
Thus, we take $\psi_0 = \psi_{d,\lambda_0}$ to be the solution of
\begin{equation} \label{e.psiddef}
\begin{aligned}
-\Delta \psi_{d,\lambda_0} &= \lambda_0 \psi_{d,\lambda_0} \quad \mbox{in $D$}\\
\psi_{d,\lambda_0} &= 1 \quad \quad \quad \mbox{at $\partial D$}.
\end{aligned}
\end{equation}
We assume here that $\lambda_0 \neq 0$ is real, and that it is not an eigenvalue of $-\Delta$ in $D$ with Dirichlet
boundary condition $0$. Under these conditions the solution of \eqref{e.psiddef} exists and is unique and real-valued.
Since we have only assumed that $D$ is a Lipschitz domain, $\psi_{d,\lambda_0}$ is in $H^2_{loc}(D) \cap H^1(D)$, which
is enough for our purposes. (In \cite{KV1} the subscript $d$ stood for ``dopant;'' here it is just a reminder that
$\psi_{d,\lambda_0}$ depends on both $D$ and $\lambda_0$.)

We note for future reference that with the substitution $\psi_0 = \psi_{d,\lambda_0}$, our normalization
\eqref{normalization-using-psi0} has become
\begin{equation} \label{e.normalize1}
\int_\Omega u_\de (x)u_0(x) \, dx = |\ENZ| + \int_D \psi_{d,\lambda_0}^2\,dx \, .
\end{equation}

Since the eigenvalues of an elliptic operator are discrete, we expect that only certain choices of
$\lambda_0$ should be acceptable. When we consider the expansion term-by-term in Section \ref{subsec:leading-order-terms},
the condition on $\lambda_0$ will emerge as the consistency condition for the existence of $\phi_1$; therefore we
like to call it the \emph{consistency condition}. However the same condition can be derived as follows:
it is easy to see from \eqref{e.evp} that
\begin{equation} \label{e.meanzero}
\int_\Omega u_\de\,dx  = 0
\end{equation}
by integrating the PDE over $\Omega$ and using the homogeneous Neumann boundary condition (along with the assumption
that $\lambda_0 \neq 0$). At leading order this gives
\begin{equation} \label{consistency-first-version}
|\ENZ| + \int_D \psi_{d,\lambda_0} \, dx = 0 \, .
\end{equation}
We shall discuss the solvability of this condition in Section \ref{subsubsec:consistency}, but we note here that
(i) it requires $\int_D \psi_{d,\lambda_0} \, dx$ to be negative, and (ii) when this integral is negative,
\eqref{consistency-first-version} simply determines the area of the ENZ region.

We assumed above that $\lambda_0$ is real-valued and nonzero. Actually, the consistency condition
\eqref{consistency-first-version} implies that it must be positive. Indeed, using the definition
\eqref{e.psiddef} of $\psi_{d,\lambda_0}$ we have
\begin{equation*}
\begin{aligned}
\int_D |\nabla \psi_{d,\lambda_0} |^2 \, dx &= \int_D \mathrm{div}\, (\psi_{d,\lambda_0} \nabla  \psi_{d,\lambda_0} ) \, dx -
\int_D\psi_{d,\lambda_0} \Delta \psi_{d,\lambda_0}\, dx\\
&= \int_{\partial D} \partial_{\nu_D} \psi_{d,\lambda_0} \, d\sh^1 + \lambda_0 \int_D \psi_{d,\lambda_0}^2 \, dx \, .
\end{aligned}
\end{equation*}
But using the PDE again along with the consistency condition we have
$$
\int_{\partial D} \partial_{\nu_D} \psi_{d,\lambda_0} \, d\sh^1 = \int_D \Delta \psi_{d,\lambda_0} \, dx
= -\lambda_0 \int_D \psi_{d,\lambda_0} \, dx = \lambda_0 |\ENZ| \, .
$$
Combining these relations, we conclude that
$$
\int_D |\nabla \psi_{d,\lambda_0} |^2 \, dx = \lambda_0 \left(|\ENZ| + \int_D \psi_{d,\lambda_0}^2 \, dx \right) \, .
$$
So $\lambda_0$ must be positive, as asserted.

\subsubsection{Statement of our theorem on analyticity of $u_\delta$ and $\lambda_\delta$} \label{subsubsection:theorem-stmt}
We are ready to state our result on the existence of eigenvalues and eigenvectors of \eqref{e.evp} depending analytically
on $\delta$.

\begin{theorem} \label{t.main1}
Let $D$ and $\Omega$ be as discussed at the beginning of Section \ref{subsec:preliminaries-etc} and let $\lambda_0$
be a positive real number which (i) is not a Dirichlet eigenvalue of $-\Delta$ in $D$ and (ii) satisfies the
consistency condition \eqref{consistency-first-version}. Then for all complex
$\delta$ in a neighborhood of $0$ there exists a simple eigenfunction $u_\delta$ of \eqref{e.evp} with eigenvalue
$\lambda_\delta$ such that
$$
\lim_{\delta \rightarrow 0} \lambda_\delta = \lambda_0 \quad \mbox{and} \quad
\lim_{\delta \rightarrow 0} u_\delta = u_0 =
\left\{
\begin{array}{cc}
1 & x \in \ENZ \\
\psi_{d,\lambda_0} & x \in D \, .
\end{array}
\right.
$$
Moreover, with the normalization \eqref{e.normalize1}
the eigenfunction $u_\delta$ and its eigenvalue $\lambda_\delta$ are complex analytic functions of $\delta$
in a neighborhood of $\delta = 0$.
\end{theorem}

To be sure the final statement is clear: we will show, in the course of the proof, that the map
$\delta \mapsto u_\delta$ is a complex analytic function of $\delta$ (near $\delta = 0$) taking values
in a suitable Banach space. This is equivalent to the statement that $u_\delta$ has a Taylor expansion
$u_\delta = u_0 + \delta u_1 + \delta^2 u_2 + \cdots$ with a positive radius of convergence (see e.g.
\cite{D} or \cite{Whittlesey}).

\subsubsection{On satisfying the consistency condition} \label{subsubsec:consistency}

Our consistency condition \eqref{consistency-first-version} involves the function $\psi_{d,\lambda_0}$, so its dependence on
$D$ is not very explicit. This subsection examines how it can be satisfied, either (a) by choosing $D$ and $\Omega$ appropriately
with $\lambda_0$ held fixed, or (b) by choosing $\lambda_0$ appropriately, with $D$ and $\Omega$ held fixed. (This discussion is not
used in the proof of Theorem \ref{t.main1}. A reader who is mainly interested in that theorem can skip to Section
\ref{subsec:leading-order-terms}.)

We start with a representation formula for $\psi_{d,\lambda_0}$ in terms of the Dirichlet eigenvalues and
eigenfunctions of the domain $D$ (more precisely, the spectrum of $-\Delta$ in $H^1_0(D)$). Let $\{\mu_n\}_{n=1}^\infty$ be
the Dirichlet eigenvalues of $-\Delta$, and let $\{\chi_n\}_{n=1}^\infty$ be an associated set of orthonormal eigenfunctions;
as usual, the eigenvalues are enumerated in nondecreasing order and repeated according to multiplicity. By expressing
the function $\psi_{d,\lambda_0} - 1 $ (which vanishes at $\partial D$) in terms of the eigenfunction basis, one finds by a routine
calculation that
\begin{equation} \label{e.psidformula}
\psi_{d,\lambda_0} = 1 + \lambda_0 \biggl( \sum_{n=1}^\infty \frac{\int_D \chi_n\,dx}{\mu_n - \lambda_0}\chi_n\biggr) \, .
\end{equation}

We see from this formula that $\psi_{d,\lambda_0}$ depends only on the eigenfunctions for which $\int_D \chi_n \, dx \neq 0$.
The following simple proposition assures us that there are infinitely many of these. (For a more quantitative result -- estimating
how many of the first $n$ Dirichlet eigenfunctions have nonzero mean -- see \cite{steinerberger-venkatraman}.)

\begin{proposition} \label{p.count}
For any bounded domain $D\subset \RR^d$ with Lipschitz boundary, there are infinitely many Dirichlet eigenfunctions
$\chi_n$ such that $\int_D \chi_n(x) \, dx \neq 0$.
\end{proposition}

\begin{proof}
The functions $\{\chi_n\}$ form an orthonormal basis of $L^2(D)$. The constant function $1$ is in $L^2(D)$ (since $D$ is bounded),
so
$$
1 = \sum_{n=1}^\infty \biggl(\int_D \chi_n(x)\,dx\biggr)\chi_n \, ,
$$
where the series on the right, if infinitely many of terms are nonzero, must be understood in the sense of convergence of
$L^2$ functions. Now, if all but finitely many of the coefficients $\int_D \chi_n(x) \, dx$ were to vanish then the constant
function $1$ would be a finite sum of eigenfunctions that all vanish at $\partial D$. This is not possible, so the
proposition is proved.
\end{proof}

The consistency condition \eqref{consistency-first-version} involves just the \emph{integral} of $\psi_{d,\lambda_0}$, which
by \eqref{e.psidformula} is
\begin{equation} \label{integral-of-psi-repn}
\int_D \psi_{d,\lambda_0} = |D| + \lambda_0 \sum_{n=1}^\infty \frac{(\int_D \chi_n(x)\,dx)^2}{\mu_n - \lambda_0} \, .
\end{equation}
We note that the sum on the right hand side of \eqref{integral-of-psi-repn} is absolutely convergent. Indeed, each term is
finite (since by hypothesis $\lambda_0$ is not a Dirichlet eigenvalue), and for all but finitely many terms $\mu_n > \lambda_0$
(since the eigenvalues are ordered and tend to infinity). Since all but finitely many of the terms are positive, the sum
converges absolutely.

We turn now to the question how the consistency condition \eqref{consistency-first-version} can be satisfied by
choosing $D$ and $\Omega$ appropriately, for any fixed positive $\lambda_0$. As already noted earlier, we need
only ask how $\int_D \psi_{d,\lambda_0} \ dx$ can be made negative, since the consistency condition is then satisfied by
choosing the area of $\ENZ$ correctly. For a given domain $D$, it is of course possible for $\int_D \psi_{d,\lambda_0} \ dx$
to be positive. However, if the Dirichlet eigenvalues of $-\Delta$ in $D$ are $\{ \mu_n \}$, then the Dirichlet eigenvalues
of $-\Delta$ in the scaled domain $tD$ are $\mu_n/t^2$. As $t$ varies, there will be selected values where
$\mu_n/t^2$ crosses $\lambda_0$ for some eigenvalue $\mu_n$ such that $\int \chi_n \, dx \neq 0$. As this crossing happens, we
see from \eqref{integral-of-psi-repn} that the value of $\int_D \psi_{d,\lambda_0} \, dx$ jumps from $-\infty$ (as $\mu_n/t^2$
approaches $\lambda_0$ from below) to $+ \infty$ (as $\mu_n/t^2$ increases past $\lambda_0$). As $t$ ranges over the interval
between two consecutive crossings, $\int_D \psi_{d,\lambda_0} \, dx$ takes every real value by the
intermediate value theorem. Thus, the scale factor $t$ can easily be chosen so that $\int_D \psi_{d,\lambda_0} \ dx$ is negative.

Finally, we examine how the consistency condition can be satisfied by choosing $\lambda_0$ appropriately when $D$ and $\Omega$
are held fixed. Combining \eqref{consistency-first-version} with \eqref{integral-of-psi-repn}, this amounts to studying the roots
of $|\Omega| + f(t) = 0 $, where
\begin{equation} \label{defn-of-f}
f(t) := t \sum_{n=1}^\infty \frac{\bigl( \int_D \chi_{n})^2}{\mu_n - t} \, .
\end{equation}
We noted above that this sum converges absolutely provided that $t$ is not an eigenvalue with a nonzero-mean eigenfunction.
Differentiating term-by-term gives
\begin{equation} \label{deriv-of-f}
f'(t) = \sum_{n=1}^\infty \frac{\mu_n \bigl( \int_D \chi_{n})^2}{(\mu_n - t)^2} > 0 \, .
\end{equation}
(This calculation is legitimate, since the differentiated sum again converges absolutely; indeed, for large $n$
its $n$th term is comparable to that of $f$.) Remembering that an eigenvalue $\mu_n$ participates in these sums only
if $\int_D \chi_n \, dx \neq 0$, it is convenient to let $J = \{ m_1, m_2, \ldots \}$ be the ordered list of eigenvalues
having at least one eigenfunction with nonzero mean (which is infinite, by Proposition \ref{p.count}).
Then we see from \eqref{defn-of-f} -- \eqref{deriv-of-f} that $f(t)$ increases monotonically from $-\infty$ to $+\infty$
on each interval $m_i < t < m_{i+1}$. Thus: each of these intervals contains a unique choice of $\lambda_0$ for which
the consistency condition holds.

\subsection{The leading order terms} \label{subsec:leading-order-terms}

We have as yet determined only the zeroth-order terms in the Taylor expansions of $\lambda_\delta$ and
$u_\delta$. We turn now to the identification of additional terms. The
first few, which are discussed in this section, are used in our proof of Theorem \ref{t.main1}; briefly, knowing them
lets us desingularize the PDE problem, permitting application of the implicit function theorem.
The higher-order terms, which we discuss in Section \ref{subsec:higher-order-terms}, are also interesting.
Indeed, the process by which the expansion is determined term-by-term is intimately related to our
implicit-function-theorem-based proof of Theorem \ref{t.main1}: the implicit function theorem requires
the invertibility of a certain linear operator, whereas our identification of each successive term in
the expansion involves inverting this operator. (We note in passing that the higher-order terms can also
be used to provide an alternative proof of Theorem \ref{t.main1} by directly majorizing the
Taylor expansions. For arguments of this type in closely analogous settings see
\cite{Fortes-Lipton-Shipman-2010,Fortes-Lipton-Shipman-2011}.)

Before delving into the details, let us provide a big-picture view of the calculation. Our plan is to substitute the
expansions \eqref{e.lambdaexp} and \eqref{e.ansatz} into the PDE \eqref{e.evp} and the normalization \eqref{e.normalize1}
and expand in powers of $\delta$. The condition that $\lambda_0$ must satisfy -- \eqref{consistency-first-version} -- will
emerge naturally as the consistency condition for the PDE problem (in $\ENZ$) that determines $\phi_1$. When this
consistency condition holds, $\phi_1$ is determined only up to an additive constant, which we call $e_1$. The function
$\psi_1$ solves a different PDE problem (in $D$), which involves $\phi_1$ and $\lambda_1$; as a result,
$\psi_1$ is initially found in terms of the not-yet-determined parameters $e_1$ and $\lambda_1$.
Finally, $e_1$ and $\lambda_1$ are determined by the normalization condition \eqref{e.normalize1} and
the consistency condition for the existence of $\phi_2$. The process by which $\phi_j$, $\psi_j$, and
$\lambda_j$ are determined for each successive $j=2,3,\ldots$ is similar.

As we shall see, each function $\phi_j$ satisfies a Poisson-type equation in $\ENZ$. The associated
consistency condition comes from the fact that if $\Delta \phi = f$ in a domain and $\partial _\nu \phi = g$
at its boundary, then the volume integral of $f$ must equal the boundary integral of $g$. When the equation
is consistent, the solution is determined only up to an additive constant. For this reason, it will be
convenient to view each $\phi_j$ as the sum of a mean-zero function and a constant:

\begin{align} \label{e.phindecomp}
	\phi_n = \ringphi_n + e_n, \quad e_n \in \RR, \quad \int_{\ENZ} \ringphi_n\,dx = 0\,.
\end{align}
This decomposition induces one of $\psi_j$, since (as we'll see) $\psi_j$ solves a linear PDE in $D$ with
$\psi_j = \phi_j$ at $\partial D$. While the form of this decomposition will emerge naturally later, we
mention it now as a complement to \eqref{e.phindecomp}:
\begin{align} \label{e.psindecomp}
	\psi_n = \ringpsi_n + e_n \psi_{d,\lambda_0}\, \quad \ringpsi_n = \ringphi_n \quad \mbox{ on } \partial D\, .
\end{align}

We turn now to the business of this subsection: identification of the initial terms in the expansions. We have, of course,
already chosen the order-one terms in the expansion of $u_\delta$:
\begin{align*}
\phi_0(x) &:= 1 \quad \mbox{ for } x \in \Omega \setminus \overline{D}, \\
\psi_0(x) &:= \psi_{d,\lambda_0}(x) \quad \mbox{ for } x \in D\,,
\end{align*}
Plugging the expansions into the PDE \eqref{e.evp}, at order one in the ENZ region we get a PDE problem for $\phi_1$:
\begin{equation} \label{e.neumann}
\begin{aligned}
-\Delta \phi_1 &= \lambda_0 \quad \mbox{ in } \Omega \setminus \overline{D}\\
\partial_\nu \phi_1 &= 0 \quad \mbox{ on } \partial \Omega\\
\partial_\nu \phi_1 &= \partial_\nu \psi_{d,\lambda_0} \quad \mbox{ on } \partial D\,.
\end{aligned}
\end{equation}
Existence requires a consistency condition. Recalling that $\nu_D$ denotes the unit normal to $\partial D$
\textit{pointing out of $D$}, the consistency condition is
\begin{equation} \label{e.consistency0}
\lambda_0|\Omega \setminus \overline{D} | = \int_{\partial D} \frac{\partial \psi_{d,\lambda_0}}{\partial \nu_D}\,d\sh^1 \, .
\end{equation}
When this holds, the solution exists but it is unique only up to an additive constant. Therefore we decompose
$$
\phi_1 = \ringphi_1 + e_1, \quad e_1 \in \RR, \quad \int_{\ENZ} \ringphi_1\,dx = 0
$$
and recognize that while $\ringphi_1$ is uniquely determined, $e_1$ is still unknown. (We note that $\ringphi_1$
is real, since $\psi_{d,\lambda_0}$ is real; therefore our assumption that $e_1$ take real values is quite natural.)

The consistency condition \eqref{e.consistency0} is equivalent to the condition on $\lambda_0$ that we introduced earlier,
namely \eqref{consistency-first-version}. Indeed, since
$$
\int_{\partial D} \frac{\partial \psi_{d,\lambda_0}}{\partial \nu_D}\,d\sh^1  =
\int_D \Delta \psi_{d,\lambda_0} \, dx = - \lambda_0 \int_D \psi_{d,\lambda_0} \, dx \, ,
$$
\eqref{e.consistency0} can be rewritten
$$
\lambda_0 |\ENZ| = - \lambda_0 \int_D \psi_{d,\lambda_0} \, dx \, ,
$$
which is equivalent to \eqref{consistency-first-version} since we always assume $\lambda_0 \neq 0$.

We turn now to the identification of the function $\psi_1$ and the constants $\lambda_1$ and $e_1$. Since we have already
considered the order-one PDE, both in $D$ (in defining $\psi_{d,\lambda_0}$) and in $\ENZ$ (in finding $\ringphi_1$), we naturally turn
to the order-$\delta$ problem in $D$. It says
\begin{align} \label{e.psi1}
	-\Delta \psi_1 &= \lambda_1 \psi_{d, \lambda_0} + \lambda_0 \psi_1 \quad \mbox{ in } D\\ \label{e.bc}
	\psi_1 &= \phi_1 \quad \mbox{ on } \partial D \, .
\end{align}
(The boundary condition comes from the fact that $u_\delta$ cannot jump across $\partial D$.) For given $\phi_1$ and $\lambda_1$,
this boundary value problem has a unique solution (since $\lambda_0$ is not a Dirichlet eigenvalue of $-\Delta$ in $D$). Since the additive
constant $e_1$ in $\phi_1$ has not yet been determined, it is convenient to make the dependence of $\psi_1$ on $e_1$ more
explicit. Therefore we decompose $\psi_1$ as in \eqref{e.psindecomp}:
\begin{equation} \label{decomp-of-psi1}
\psi_1 = \ringpsi_1 + e_1 \psi_{d,\lambda_0} \, ,
\end{equation}
where
\begin{equation}\label{e.ringpsi1def}
\begin{aligned}
	-\Delta \ringpsi_1 &= \lambda_1 \psi_{d, \lambda_0} + \lambda_0 \ringpsi_1 \quad \mbox{ in } D\\
	\ringpsi_1 &= \ringphi_1 \quad \mbox{ on } \partial D \, .
\end{aligned}
\end{equation}
Since $\lambda_0$ is not a Dirichlet eigenvalue of $D$, any choice of $\lambda_1$ uniquely determines $\ringpsi_1$.

We must still determine $e_1$ and $\lambda_1$. For this purpose, we shall use the normalization condition \eqref{e.normalize1} and
the condition that
\begin{equation} \label{e.ringpsi1mean}
\int_D \ringpsi_1\,dx  = 0 \, .
\end{equation}
Some explanation is in order about the latter. Remember that while our condition on $\lambda_0$ was initially obtained by
requiring that $\int_\Omega u_\delta \, dx = 0$ at order one, the same condition emerged above as the
consistency condition for existence of $\phi_1$. The status of \eqref{e.ringpsi1mean} is similar. It is at once

\begin{enumerate}
\item[(a)] the order-$\delta$ version of the condition that $\int_\Omega u_\delta \, dx = 0$, and

\item[(b)] the consistency condition for existence of $\phi_2$.
\end{enumerate}
To see (a), we observe that
\begin{align*}
\int_{\ENZ} \phi_1 \, dx + \int_D \psi_1 \, dx &= e_1 |\ENZ| + \int_D \left( \ringpsi_1 + e_1 \psi_{d,\lambda_0} \right) \, dx \\
&= \int_D \ringpsi_1 \, dx
\end{align*}
using the consistency condition \eqref{consistency-first-version} in the second line. We postpone the justification of (b) to the end of
this subsection, since it requires a bit of calculation and it isn't immediately needed.

We now identify the value of $\lambda_1$. Multiplying both sides of \eqref{e.ringpsi1def} by $\psi_{d,\lambda_0}$ and
integrating gives
\begin{equation} \label{e.derivation_of_lambda1}
\begin{aligned}
   &\lambda_1 \int_D  \psi_{d,\la_0}(\psi_{d,\la_0} - 1)\,dx =
    - \int_D \bigl(\Delta \ringpsi_1 + \la_0 \ringpsi_1\bigr)(\psi_{d,\la_0} - 1) \,dx\\
   &\quad \quad  \underset{\eqref{e.psiddef},\eqref{e.ringpsi1mean}}{=}
     \int_{\partial D} \ringphi_1 \partial_{\nu_D}\psi_{d,\la_0}\,d\sh^1\\
   &\quad \qquad \quad \underset{\eqref{e.neumann}}{=} \int_{\partial D} \ringphi_1 \partial_{\nu_D}\ringphi_1\,d\sh^1\\
   &\quad \qquad \quad \underset{\eqref{e.neumann}}{=} -\int_{\ENZ} |\nabla \ringphi_1|^2 \, dx \, .
\end{aligned}
\end{equation}
Combining this with \eqref{consistency-first-version}, we conclude that
\begin{equation} \label{e.lambda1def}
\lambda_1 = - \frac{1}{ |\ENZ| + \int_D \psi_{d,\la_0}^2 } \int_{\ENZ}|\nabla \ringphi_1|^2 \, dx \, .
\end{equation}
We note that this $\lambda_1$ does not depend on the as-yet undetermined constant $e_1$.

Finally, we identify the value of $e_1$ using the order $\delta$ term in the expansion of the normalization condition
\eqref{e.normalize1}, which is
$$
\int_D \psi_1 \psi_{d,\lambda_0} \, dx + \int_{\ENZ} \phi_1 = 0.
$$
Using \eqref{decomp-of-psi1}, this is equivalent to
$$
\int_D \ringpsi_1 \psi_{d,\lambda_0} + e_1 \psi_{d,\lambda_0}^2\, dx + e_1 |\ENZ| = 0 \, ,
$$
so
\begin{equation} \label{e.e1def}
e_1 = - \frac{\int_D \ringpsi_1 \psi_{d,\lambda_0}\,dx}{|\ENZ| + \int_D \psi_{d,\lambda_0}^2\,dx} \, .
\end{equation}
We note that this definition is not circular: the right hand side of \eqref{e.e1def} involves $\ringpsi_1$, which is defined by
\eqref{e.ringpsi1def} and which therefore depends on $\lambda_1$. However $\ringpsi_1$ is independent of $e_1$, since our chosen
value of $\lambda_1$ -- given by \eqref{e.lambda1def} -- is independent of $e_1$.

A thoughtful reader might ask: is it really true that $\int_D \ringpsi_1 \, dx = 0$ when $\lambda_1$ is given by
\eqref{e.lambda1def} and $\ringpsi_1$ is determined by \eqref{e.ringpsi1def}? The answer is yes. To see why, we revisit the
the calculation \eqref{e.derivation_of_lambda1} without assuming that this integral vanishes:
\begin{equation} \label{e.meanzerorigorous}
\begin{aligned}
 &\lambda_1 \int_D  \psi_{d,\la_0}(\psi_{d,\la_0} - 1)\,dx =
   - \int_D \bigl(\Delta \ringpsi_1 + \la_0 \ringpsi_1\bigr)(\psi_{d,\la_0} - 1) \,dx\\
 &\quad \quad  \underset{\eqref{e.ringpsi1def},\eqref{e.psiddef}}{=}
   \int_{\partial D} \ringphi_1 \partial_{\nu_D}\psi_{d,\la_0}\,d\sh^1 + \lambda_0 \int_D \ringpsi_1 \,dx \\
 &\quad \qquad \quad \underset{\eqref{e.neumann}}{=} -\int_{\ENZ} |\nabla \ringphi_1|^2\,dx + \lambda_0 \int_D \ringpsi_1 \,dx \, .
\end{aligned}
\end{equation}
This amounts to a linear relation between $\lambda_1$ and $\int_D \ringpsi_1 \, dx$ (since $\lambda_0$, $\psi_{d,\lambda_0}$, and
$\ringphi_1$ are by now fixed). Our choice of $\lambda_1$ is precisely the one that makes $\int_D \ringpsi_1 \, dx$ vanish.

We note for future reference that the functions $\phi_1$ and $\psi_1$ satisfy the order-$\delta$ versions of
$\int_\Omega u_\delta \, dx = 0$ and our normalization condition \eqref{e.normalize1}, namely
\begin{equation} \label{e.consistency-base} \int_{\ENZ} \phi_1 + \int_D\psi_1 \,dx = 0\,.
\end{equation}
and
\begin{equation} \label{e.norm-basecase}
\int_{\ENZ} \phi_1 \,dx + \int_D\psi_1 \psi_{d,\la_0} = 0\,.
\end{equation}
(Indeed, we found $\lambda_1$ and $e_1$ by assuring these relations.)

We close this subsection by justifying our claim that the condition $\int_D \ringpsi_1 \, dx = 0$ is equivalent to the consistency
condition for existence of $\phi_2$. Our starting point is the PDE for $\phi_2$,
which is the order-$\delta$ PDE in $\ENZ$:
$$
\begin{aligned}
-\Delta \ringphi_{2} &=  \lambda_{0} \phi_{1} + \lambda_1 \phi_0\, \quad &\mbox{ in } \Omega \setminus \overline{D}\\
\partial_{\nu_\Omega}\ringphi_{2} &= 0  \quad &\mbox{ on } \partial \Omega\\
\partial_{\nu_D} \ringphi_{2} &= \partial_{\nu_D} \psi_1 =
\partial_{\nu_D} \ringpsi_1 + e_1 \partial_{\nu_D} \psi_{d,\lambda_0}\quad &\mbox{ on } \partial D \, .
\end{aligned}
$$
Its consistency condition (remembering that $\nu_D$ points outward from $D$) is
$-\int_{\ENZ} \Delta \phi_2 \, dx = \int_{\partial D} \partial_{\nu_D} \phi_2 \, d\sh^1 $, in other words
\begin{equation} \label{consistency-for-phi2}
\lambda_0 e_1 |\ENZ| + \lambda_1 |\ENZ| =
\int_{\partial D} \partial_{\nu_D} \ringpsi_1 + e_1 \partial_{\nu_D} \psi_{d,\lambda_0} \, d\sh^1 \, .
\end{equation}
The right side is equal to
$$
\int_D \Delta \ringpsi_1 + e_1 \Delta \psi_{d, \lambda_0} \, dx =
- \int_D \left (\lambda_1 \psi_{d,\lambda_0} + \lambda_0 \ringpsi_1 \right) \, dx - e_1 \int_D \lambda_0 \psi_{d,\lambda_0} \, dx \, .
$$
Using this along with the consistency condition ($|\ENZ| + \int_D \psi_{d,\lambda_0} \, dx = 0$), \eqref{consistency-for-phi2}
reduces to
$$
\lambda_0 \int_D \ringpsi_1 \, dx = 0 \, ,
$$
which demonstrates our claim (since $\lambda_0 \neq 0$).

\subsection{The higher order terms} \label{subsec:higher-order-terms}
In this section we explain how the remaining terms in the expansions for $u_\delta$ and $\lambda_\delta$ can be found by
an inductive procedure. As noted earlier, it is possible to prove Theorem~\ref{t.main1} by majorizing the resulting series.
However our proof -- presented in Section \ref{subsec:analyticity} --
uses a different approach, based on the implicit function theorem.
Therefore the material in this section will not be used in the rest of the paper;
a reader who is mainly interested
in the proof of Theorem \ref{t.main1} can skip directly to Section \ref{subsec:analyticity}.

Our procedure is inductive: given
\begin{align}
\label{e.ind.hyp}
\{\ringphi_j\}_{j=1}^n \, \, ,  \{\lambda_j\}_{j=1}^{n} \, \, , \{\ringpsi_j\}_{j=1}^{n}\, \, , \{e_j\}_{j=1}^{n},
\end{align}
satisfying certain properties, we shall explain how to find $\ringphi_{n+1}$, $\lambda_{n+1}$, $\ringpsi_{n+1}$,
and $e_{n+1}$ with the analogous properties at level $n+1$. The base case of the induction will be
$j=1$, which was addressed in the previous subsection. Throughout this discussion, we understand that
$\phi_j$ and $\psi_j$ are determined by $\ringphi_j$, $\lambda_j$, $\ringpsi_j$, and $e_j$ via
\begin{equation} \label{e.decomp}
\phi_j := \ringphi_j + e_j, \quad \psi_j = \ringpsi_j + e_j \psi_{d,\lambda_0} \, .
\end{equation}
\medskip

\noindent \emph{Inductive hypotheses:}
\begin{itemize}

\item \emph{For $j = 1, \ldots, n$, the functions $\ringphi_j$ and $\ringpsi_j$ satisfy}
\begin{equation} \label{e.psiringj0}
\int_{\ENZ} \ringphi_j \, dx = 0 \quad \mbox{\emph{and}} \quad \int_D \ringpsi_j \, dx = 0 \, .
\end{equation}
We note that when $\ringphi_j$ has mean zero and $\lambda_0$ satisfies the consistency
condition \eqref{consistency-first-version},
$$
\int_{\ENZ} \phi_j \,dx + \int_D \psi_j \, dx = \int_D \ringpsi_j \, dx \, ;
$$
thus, the condition that $\ringpsi_j$ have mean zero is equivalent to
\begin{equation} \label{e.normconseq2}
\int_{\ENZ} \phi_j \,dx + \int_D \psi_j \, dx = 0 \, ,
\end{equation}
which amounts to the condition that $\int_\Omega u_\delta \, dx = 0$ at order $\delta^j$.

\item \emph{For $j = 1, \ldots, n$ the constant $e_j$ is chosen so that}
\begin{equation} \label{e.normconseq1}
\int_{\ENZ} \phi_j \,dx + \int_D \psi_j\psi_{d,\lambda_0}\,dx = 0 \, .
\end{equation}
This is simply our normalization condition \eqref{e.normalize1} at order $\delta^j$. Using
\eqref{e.decomp}, we see that it is equivalent to
\begin{equation} \label{choice-of-ej}
e_j = - \frac{1}{|\ENZ| + \int_D \psi_{d,\lambda_0}^2\,dx} \int_{D}\ringpsi_j \psi_{d,\lambda_0}\,dx \, .
\end{equation}
\end{itemize}

A useful identity arises by combining \eqref{e.normconseq2} and \eqref{e.normconseq1}: subtracting one from the other
gives the orthogonality relation
\begin{equation} \label{e.cancel}
\int_D \psi_j (1 - \psi_{d,\lambda_0}) \,dx = 0\, .
\end{equation}

Of course, the functions $\ringphi_j$, $\ringpsi_j$ and the constants $e_j$ and $\lambda_j$ will be chosen for
$j = 1, \ldots, n$ so that the associated expansions satisfy our PDE \eqref{e.evp} to a certain order, and
the inductive step (choosing these quantities for $j=n+1$) will assure that the PDE is satisfied to the next order.

As already noted, the base case $j=1$ is already in place. Indeed, the functions $\ringphi_1$, $\ringpsi_1$ and the
constants $e_1, \lambda_1$ found in Section \ref{subsec:leading-order-terms} have the desired properties (see
\eqref{e.norm-basecase} and \eqref{e.consistency-base}).
\medskip

\noindent \emph{Induction step:} We will determine $\ringphi_{n+1}$ using the PDE in the ENZ region
at order $\delta^n$; then we will determine $\ringpsi_{n+1}$, $e_{n+1}$, and $\lambda_{n+1}$ by using
the PDE in the region $D$ at order $\delta^{n+1}$ combined with
the $j=n+1$ versions of conditions \eqref{e.psiringj0} and \eqref{e.normconseq1}.

Since our argument uses the entire expansion of $u_\delta$ and $\lambda_\delta$, we take the convention that
$\ringphi_0 = 0$, $e_0 = 1$ so that $\phi_0 = \ringphi_0 + e_0  = 1$; similarly, we take $\ringpsi_0 = 0$ so that
$\psi_0 = \psi_{d,\lambda_0}$.

The function $\ringphi_{n+1} \in H^1(\ENZ)$ is obtained by substituting the expansion into the PDE, then focusing on the equation
in $\ENZ$ at order $\delta^n$. This gives the Neumann problem
\begin{equation}\label{e.phin+1}
\begin{aligned}
-\Delta \ringphi_{n+1} &= \sum_{k=0}^{n} \lambda_{k} \phi_{n-k}\, \quad &\mbox{ in } \Omega \setminus \overline{D}\\
\partial_{\nu_\Omega}\ringphi_{n+1} &= 0  \quad &\mbox{ on } \partial \Omega\\
\partial_{\nu_D} \ringphi_{n+1} &= \partial_{\nu_D} \psi_n \quad &\mbox{ on } \partial D\\
\int_{\Omega \setminus \overline{D}} \ringphi_{n+1} (x)\,dx &= 0\, .
\end{aligned}
\end{equation}
For a solution to exist, the integral over $\ENZ$ of the bulk source term must be consistent with the integral over $\partial D$ of
$\partial_{\nu_D}\psi_n.$ Using the PDE for $\psi_n$, the boundary integral can be expressed as
a bulk integral over $D$. This leads to the consistency condition
\begin{equation} \label{e.phin+1.consistency}
 \sum_{k=0}^n \lambda_k \biggl[ \int_{\ENZ} \phi_{n-k}\,dx +  \int_D \psi_{n-k}\,dx \biggr]= 0 \, ,
\end{equation}
which holds thanks to \eqref{e.normconseq2}. Thus the PDE problem \eqref{e.phin+1} is consistent, and
$\ringphi_{n+1}$ is its unique mean-zero solution.

Turning now to the PDE in $D$, at order $\delta^{n+1}$ we find the Dirichlet problem
\begin{equation}\label{e.psin+1}
\begin{aligned}
- \Delta \psi_{n+1} - \lambda_0 \psi_{n+1} &= \sum_{k=1}^{n+1} \lambda_k \psi_{n+1-k} \quad & \mbox{ in } D\\
\psi_{n+1} &= \phi_{n+1} = \ringphi_{n+1} + e_{n+1}\psi_{d,\lambda_0}\quad & \mbox{ on } \partial D\,.
\end{aligned}
\end{equation}
By linearity (and using the definition of $\psi_{d,\lambda_0}$), it suffices to solve
\begin{equation} \label{e.ringpsin}
\begin{aligned}
- \Delta \ringpsi_{n+1} - \lambda_0 \ringpsi_{n+1} &= \sum_{k=1}^{n+1} \lambda_k \psi_{n+1-k} \quad & \mbox{ in } D\\
\ringpsi_{n+1} &= \ringphi_{n+1} \quad & \mbox{ on } \partial D\,.
\end{aligned}
\end{equation}
The solution $\ringpsi_{n+1}$ depends on $\lambda_{n+1}$, which is as yet unknown. Its value can be
obtained by arguing as we did for $j=1$ in Section \ref{subsec:leading-order-terms}. Inspired by that calculation,
we \emph{define} $\lambda_{n+1}$ by
\begin{equation}
\label{e.lambdan+1def}
\lambda_{n+1} :=
\frac{1}{|\ENZ| + \int_D\psi_{d,\la_0}^2} \int_{\partial D} \ringphi_{n+1}\partial_{\nu_D} \psi_{d,\lambda_0} \, d\sh^1 \, ,
\end{equation}
which is well-defined, since the right hand side involves only quantities that have already been determined.
Since we have now fixed $\lambda_{n+1},$ the PDE \eqref{e.ringpsin} determines $\ringpsi_{n+1}.$ We may then choose $e_{n+1}$ by
\begin{equation} \label{e.en+1}
e_{n+1} := - \frac{1}{|\ENZ| + \int_D \psi_{d,\lambda_0}^2\,dx} \int_D \ringpsi_{n+1}\psi_{d,\lambda_0} \, dx \, .
\end{equation}

To complete the induction we must check that our choices for $j=n+1$ meet the requirements of the
inductive hypothesis. To do so, it suffices to check that \eqref{e.psiringj0} and \eqref{e.normconseq1} hold for $j=n+1$.
The latter follows immediately from our choice of $e_{n+1}$. To get the former, we multiply both sides of
\eqref{e.ringpsin} by $(\psi_{d,\la_0} - 1),$ remembering that this function vanishes at $\partial D$.
Integrating, using the orthogonality in \eqref{e.cancel}, and remembering our convention that
$\psi_0 = \psi_{d,\la_0},$ this calculation gives
\begin{equation*}
\begin{aligned}
\lambda_{n+1} \int_D \psi_{d,\la_0}(\psi_{d,\la_0} - 1)\,dx &= - \int_D (\Delta \ringpsi_{n+1}
 + \la_0 \ringpsi_{n+1})(\psi_{d,\la_0} - 1)\,dx \\
 &= \int_{\partial D} \ringpsi_{n+1} \partial_{\nu_D}\psi_{d,\la_0}\,d\sh^1 + \la_0 \int_D\ringpsi_{n+1}\\
 &= \int_{\partial D} \ringphi_{n+1}\partial_{\nu_D}\psi_{d,\lambda_0} \, d\sh^1  + \la_0 \int_D\ringpsi_{n+1}\, .
\end{aligned}
\end{equation*}
Combining this with the definition \eqref{e.lambdan+1def} of $\lambda_{n+1}$, and remembering that
$\int_D \psi_{d,\la_0}(\psi_{d,\la_0} - 1)\,dx = \int_D \psi_{d,\la}^2\,dx + |\ENZ|$, we conclude that
\begin{equation*}
	\int_D \ringpsi_{n+1} = 0 \, .
\end{equation*}
Since $\ringphi_{n+1}$ was chosen from the start to have mean value zero, we have confirmed \eqref{e.psiringj0} and the
induction is complete.

\subsection{The proof of Theorem~\ref{t.main1}}
\label{subsec:analyticity}

Our proof will be based on the implicit function theorem. To get started, we shall ``desingularize''
our eigenvalue problem \eqref{e.evp}, reformulating it in a way that doesn't involve dividing by
$\delta$. This will be done by using the leading order terms of the expansion.

Since we have assumed very little regularity for $\partial D$ and $\partial \Omega$ -- they are merely Lipschitz domains -- we cannot
expect the second derivatives of $u_\delta$ to be in $L^2(\Omega)$. Therefore we must work with a fairly weak solution of the PDE.
However, standard elliptic regularity results show that our $u_\delta$ is actually a smooth function of $x$ away from the boundaries
$\partial \Omega$ and $\partial D$. If the boundaries are smooth, then $u_\delta$ is also smooth up to the boundaries (though it cannot
be smooth across $\partial D$, since $\eps_\delta$ jumps there).

We now establish some notation and discuss the functional analytic framework we will use. The reader can refer to \cite{McL} for proofs
of the following facts.

\begin{itemize}
\item Given a bounded region $A \subset \RR^2$ with Lipschitz continuous boundary (in our setting, $A$ will either be $\ENZ$ or $D$),
let $\nu_A$ denote the unit normal vector field that points \emph{out} of $A$. This normal vector exists at $\mathcal{H}^1-$ almost every
point of the boundary $\partial A$, by Rademacher's theorem.

\item As usual, $H^1(A)$ denotes (possibly complex-valued) square-integrable functions on $A$ with distributional gradients that are
also represented by integration against an $L^2$ vector field. Functions in $H^1(A)$ have a boundary trace.
More precisely, there is a bounded linear operator
$\gamma_0: H^1(A) \to H^{1/2}(\partial A)$ that is surjective. It has the property that $\gamma_0(f)(x) = f(x)$ for any function
$f \in H^1 (A) \cap C(\overline{A})$, at $\sh^1-$almost every $x \in \partial A.$  When we want to indicate the dependence of
$\gamma_0$ on the domain $A$, we will write $\gamma_{0,A}$.

\end{itemize}

\noindent We shall also use the fact that if $\xi$ is an $L^2$ vector field defined on a bounded Lipschitz domain $A$
with $\nabla \cdot \xi \in L^2(A)$, then it has a well-defined normal trace $\xi \cdot \nu_A$ in
$H^{-1/2}(\partial A)$ (the dual of $H^{1/2}(\partial A)$ using the $L^2$ inner product). It is
defined by the property that for any $u \in H^1(A)$ with $\gamma_0(u) = f$,
\begin{equation} \label{defn-xi-dot-nu}
\langle \xi \cdot \nu_A , f \rangle_{H^{-1/2}(\partial A) \times H^{1/2}(\partial A)} =
\int_A (\nabla \cdot \xi) u + \xi \cdot \nabla u \, dx \, ,
\end{equation}
and it satisfies
\begin{equation} \label{xi-dot-nu-continuous}
\| \xi \cdot \nu_A \|_{H^{-1/2}(\partial A)} \leq
C \left( \|\xi \|_{L^2(A)} + \|\nabla \cdot \xi\|_{L^2(A)} \right) \, .
\end{equation}
This is well-known, but we briefly sketch the proof since it is not very explicit in \cite{McL}.
The property \eqref{defn-xi-dot-nu} determines a well-defined linear functional on $H^{1/2}(\partial A)$
since every $f \in H^{1/2}(\partial A)$ is the boundary trace of some $u \in H^1(A)$; we use here the
fact that if $\gamma_0(u_1) = \gamma_0(u_2)$ then $u'=u_1-u_2$ can be approximated in $H^1(A)$
by compactly supported functions, so $\int_A (\nabla \cdot \xi) u' + \xi \cdot \nabla u' \, dx = 0$.
The linear functional defined this way satisfies \eqref{xi-dot-nu-continuous}, since for every
$f \in H^{1/2}(\partial A)$ there exists $u$ such that $\gamma_0(u)=f$ and
$\|u\|_{H^1(A)} \leq C \|f\|_{H^{1/2}(\partial A)}$. We will use
\eqref{defn-xi-dot-nu}--\eqref{xi-dot-nu-continuous} as follows:

\begin{itemize}

\item Let
$$
S(A) := \{ f \in H^1(A) : \Delta f \in L^2(A)\} \, ,
$$
where $\Delta f$ denotes the distributional Laplacian of $f$. Then there is a bounded linear map
(the normal derivative trace) $\gamma_1 : S(A) \to H^{-1/2}(\partial A)$. It has the property that for any $f \in C^1(\overline{A})$, $\gamma_1(f)(x) = \nu_A \cdot \nabla f(x)$ at $\sh^1-$almost every $x \in \partial A.$  The map
$\gamma_1$ is surjective, by a straightforward application of the Lax-Milgram lemma. When we want to indicate the
dependence of $\gamma_1$ on $A$, we will write $\gamma_{1,A}$.

\item There is an integration by parts formula: for any $\theta \in H^{-1/2}(\partial A)$, if $f \in S(A)$ is such that
$\gamma_1(f) = \theta$, then for all $g \in H^1(A)$ we have
$$
\langle \theta, \gamma_0(g) \rangle_{H^{-1/2}(\partial A) \times H^{1/2}(\partial A)} =
\langle \gamma_1(f), \gamma_0(g) \rangle_{H^{-1/2}(\partial A) \times H^{1/2}(\partial A)} =
\int_A (\Delta f)g + \nabla f\cdot \nabla g \, dx \, .
$$
By a convenient abuse of notation we will denote the left hand side by the more familiar expression
$\int_{\partial A}g \partial_{\nu_A} f \, d\sh^1$.

\item The following version of the divergence theorem is
obtained by taking $g = 1$ in the preceding identity:
$$
\langle \gamma_1(f), 1 \rangle_{H^{-1/2}(\partial A) \times H^{1/2}(\partial A)} = \int_A \Delta f \, dx \, .
$$
\end{itemize}

\noindent We are now ready for the proof of our main theorem.
\medskip

\begin{proof}[Proof of Theorem \ref{t.main1}.]
We break up the argument into five steps.
\medskip

\noindent {\sc Step 1:} We begin by restating our problem in a form that is amenable to use of the implicit function theorem.
To find $u_\delta$ and $\lambda_\delta$, we shall seek functions $f_\de \in H^1(\ENZ),$ and $g_\de \in H^1(D)$ and a
real number $\mu_\de$ such that
\begin{align} \label{e.proofansatz}
    u_\de := \left\{
    \begin{array}{cc}
        1  + \de f_\de &   x \in \ENZ \\
       \psi_{d,\la_0} + \de g_\de  & x \in D
    \end{array}
    \right.
\end{align}
and
\begin{align} \label{e.eigval.ansatz}
    \la_\de := \la_0 + \de \mu_\de
\end{align}
satisfy \eqref{e.evp} and the normalization \eqref{e.normalize1}. Note that in view of our formal expansion we expect
$$
f_\delta = \phi_1 + \delta \phi_2 + \ldots \, , \quad g_\delta = \psi_1 + \delta \psi_2 + \ldots \, , \quad
\mu_\delta = \lambda_1 + \delta \lambda_2 + \ldots \, ,
$$
so when $\delta = 0$ we expect
\begin{equation} \label{f0-g0-mu0}
\mbox{$f_0 = \phi_1$, $g_0 = \psi_1$, and $\mu_0 = \lambda_1$.}
\end{equation}
The point of proceeding this way is that when our PDE \eqref{e.evp} is written in terms of $f_\delta$, $g_\delta$, and $\mu_\delta$,
there are no longer any negative powers of $\delta$. (For example, the PDE $\delta^{-1} \Delta u_\delta + \lambda_\delta u_\delta = 0$
in the ENZ region $\ENZ$ becomes $\Delta f_\delta + (\lambda_0 + \delta \mu_\delta)(1 + \delta f_\delta) = 0$.)

To apply the implicit function theorem, we shall break our PDE $\nabla \cdot \eps_\delta^{-1} \nabla u_\delta + \lambda_\delta u_\delta = 0$
into three main statements: (i) the PDE holds in $D$, (ii) the PDE holds in $\ENZ$, and (iii) the continuity of
$\eps_\delta \partial_\nu u_\delta$ at $\partial D$. (There are of course other conditions:
$u_\delta$ must be continuous across $\partial D$;
$\partial_{\nu_\Omega} u_\delta$ must vanish at $\partial \Omega$; and our normalization condition must be imposed. These will be built into our
chosen function spaces.)

A first-pass idea for proceeding would be to define a function $F(\delta, \mu, f, g)$ such that our PDE is equivalent
to $F(\delta, \mu_\delta, f_\delta, g_\delta) =0$, then prove existence of $(\mu_\delta, f_\delta, g_\delta)$ by applying the
implicit function theorem. Our argument is slightly different, because we need an additional constant $c_\delta$ to satisfy a
consistency condition for the PDE in $\ENZ$. Therefore we shall
\begin{enumerate}
\item[(a)] define a function $F(\delta, \mu, f, g, c)$ such that our PDE is equivalent to
$F(\delta, \mu_\delta, f_\delta, g_\delta, 0) =0$; then we'll
\item[(b)] apply the implicit function theorem to solve $F(\delta, \mu_\delta, f_\delta, g_\delta, c_\delta) =0$; then finally
\item[(c)] we'll use the specific structure of $F$ to show that this solution has $c_\delta = 0$.
\end{enumerate}
\smallskip

\noindent {\sc Step 2:} We now make the plan concrete by specifying two Banach spaces $X$ and $Y$ and the function $F : X \rightarrow Y$
that will be used. The space $X$ is a subspace of
$$
\tilde{X}:= \CC \times \CC \times  S(\ENZ) \times S(D) \times \CC
$$
defined by
\begin{multline} \label{e.spaceX}
X := \Biggl\{ (\delta,\mu, f,g, c) \in \tilde{X} \mbox{ such that } \gamma_{0,\ENZ} (f) = \gamma_{0,D}(g) \, , \ \gamma_{1,\ENZ}(f)|_{\partial\Omega} = 0 \, ,\\
\int_{\ENZ} f + \int_D g\psi_{d,\la_0} = 0 \, \ \mbox{and} \ \int_{\ENZ} f + \int_D g = 0 \Biggr\} \, .
\end{multline}
We note that the restrictions defining $X$ assure that $u_\delta$ (determined by $f$, $g$, and $\delta$ via \eqref{e.proofansatz})
(i) does not jump across $\partial D$, (ii) satisfies our homogeneous Neumann condition at $\partial \Omega$,
(iii) satisfies our normalization condition \eqref{e.normalize1}, and (iv) satisfies $\int_\Omega u\, dx =0$.
The space $Y$ is
\begin{equation} \label{e.spaceY}
Y := L^2(\ENZ) \times H^{-1/2}(\partial D) \times L^2(D) \, .
\end{equation}
The function $F$ is defined by
\begin{equation} \label{e.functionF}
F(\de, \mu, f, g, c) := \begin{pmatrix}
        \Delta f + \lambda_0 + \delta (\mu + f(\lambda_0 + \delta \mu)) + c\\
        \partial_{\nu_D}f - \partial_{\nu_D}(\psi_{d,\la_0} + \delta g) \\
         \Delta g + \lambda_0 g + \mu (\psi_{d,\la_0} + \de g)
    \end{pmatrix} \, .
\end{equation}
Lest there be any confusion concerning the middle component: since $\nu_D$ points outward from
$D$, $\partial_{\nu_D}f$ is really $-\gamma_{1,\ENZ}(f)$. Similarly, $\partial_{\nu_D}(\psi_{d,\la_0} + \delta g)$ is really
$\gamma_{1,D}(\psi_{d,\la_0} + \delta g)$. Evidently, the middle component of $F$ is the difference between two well-defined
elements of $H^{-1/2}(\partial D)$.
\medskip

\noindent {\sc Step 3:} We will apply the implicit function theorem to get the existence of
$(\mu_\delta, f_\delta, g_\delta, c_\delta)$, depending analytically on $\delta$ near $\delta = 0$, with
$\mu_0$, $f_0$, and $g_0$ given by \eqref{f0-g0-mu0} and $c_0 = 0$. While there is a version of the
implicit function theorem in the analytic setting (see e.g. \cite{Whittlesey}), the more familiar $C^1$
version (e.g. \cite[Theorem 10.2.1]{D}) is sufficient for our purposes. Indeed, it assures the existence of
$(\mu_\delta, f_\delta, g_\delta, c_\delta)$ with continuous (complex) derivatives with respect to $\delta$.
We may then appeal to the fact that such functions are complex analytic (see e.g. \cite[Theorem 9.10.1]{D}).
It is of course crucial that
$$
F(0, \mu_0, f_0, g_0,0) =0;
$$
our choices \eqref{f0-g0-mu0} do have this property (see Sections \ref{subsec:preliminaries-etc} and \ref{subsec:leading-order-terms}).

Since our goal is to solve $F(\delta, z_\delta) = 0$ near $\delta = 0$ with $z=(\mu,f,g,c)$, we must check that (i) $F$ is
$C^1$, and that (ii) the partial differential of $F$ with respect to $z$ is invertible at $(0,z_0)$ with $z_0=(\mu_0,f_0,g_0,0)$. For (i), let us express the differential $DF$ at $(\delta, \mu, f, g, c)$ as a linear
map from $X$ to $Y$:
\begin{multline} \label{full-differential}
DF_{(\delta, \mu, f, g, c)} (\dot{\delta},\dot{\mu},\dot{f},\dot{g},\dot{c})  =
\frac{d}{dt}|_{t=0} F(\delta + t\dot{\delta}, \mu + t\dot{\mu}, f + t \dot{f}, g + t \dot{g}, c + t\dot{c} ) \\
= \begin{pmatrix}
    \dot{\delta} (\mu + f(\lambda_0 + 2 \delta \mu)) +
    \dot{\mu} (\delta + \delta^2 f) +
    \Delta \dot{f} + \dot{f} \delta (\lambda_0 + \delta \mu) + \dot{c}\\[1pt]
    -\dot{\delta} \partial_{\nu_D} g + \partial_{\nu_D} \dot{f} - \delta \partial_{\nu_D} \dot{g}\\[1pt]
    \dot{\delta} \mu g + \dot{\mu} (\psi_{d,\la_0} + \delta g) + \Delta \dot{g} + \dot{g} (\lambda_0 + \mu \delta)
    \end{pmatrix} \, .
\end{multline}
It is now straightforward to see that $DF$ depends continuously (as an operator from $X$ to $Y$) upon
$(\delta,\mu, f,g,c) \in X$. The more subtle task is point (ii). Substituting
$(\delta,\mu,f,g,c) = (0,\mu_0,g_0,g_0,0) = (0,z_0)$
in \eqref{full-differential} and taking $\dot{\delta} =0$, we see that the operator to be inverted takes
the subspace of $X$ defined by $\delta = 0$ to $Y$, mapping
$$
\dot{z} = (\dot{\mu}, \dot{f}, \dot{g}, \dot{c})
$$
to
\begin{equation} \label{e.firstvar}
D_zF_{(0,z_0)} (\dot{\mu},\dot{f},\dot{g},\dot{c}) =
\begin{pmatrix}
    \Delta \dot{f} + \dot{c}\\
    \partial_{\nu_D} \dot{f}\\
    \Delta \dot{g} + \la_0 \dot{g} + \psi_{d,\la_0} \dot{\mu}
    \end{pmatrix} \, .
\end{equation}
So our task is to prove that for all $p,q,r \in Y$, the linear system
\begin{equation} \label{linear-system}
\begin{aligned}
    \Delta \dot{f} + \dot{c} &= p \quad \mbox{ in } \ENZ\\
    \partial_{\nu_D} \dot{f} &= q \quad \mbox{ on } \partial D\\
    \Delta \dot{g} + \la_0 \dot{g} + \psi_{d,\la_0} \dot{\mu} &= r \quad \mbox{ in } D
\end{aligned}
\end{equation}
has a unique solution $(\dot{\mu}, \dot{f}, \dot{g}, \dot{c})$ in $\CC \times S(\ENZ) \times S(D) \times \CC$ satisfing
$$
\gamma_{0,\ENZ} (\dot{f}) = \gamma_{0,D}(\dot{g}) \, , \gamma_{1,\ENZ}(\dot{f})|_{\partial\Omega} = 0 \, ,
\int_{\ENZ} \dot{f} + \int_D \dot{g}\psi_{d,\la_0} = 0 \, \ \mbox{and} \ \int_{\ENZ} \dot{f} + \int_D \dot{g} = 0 \, ,
$$
and that the solution operator (the map taking $(p,q,r)$ to $(\dot{\mu},\dot{f},\dot{g},\dot{c})$) is a bounded linear map from $Y$ to
$\CC \times S(\ENZ) \times S(D) \times \CC$.

The execution of this task is, of course, very similar to the method by which we found $\lambda_2$, $\phi_2$, and $\psi_2$ in
Section \ref{subsec:higher-order-terms}. We begin by considering the first two equations in \eqref{linear-system}, which give
a PDE for $\dot{f}$ in the ENZ region $\ENZ$ with a Neumann boundary condition. For a solution to exist, the consistency condition
\begin{equation*}
- \int_{\partial D}  q \,d \sh^1 +  \dot{c} \, |\ENZ|  = \int_{\ENZ} p \, dx
\end{equation*}
must hold; therefore the solution has
\begin{equation*}
\dot{c} = \frac{1}{|\ENZ|} \biggl( \int_{\ENZ} p \, dx + \int_{\partial D} q \,d\sh^1\biggr)\,,
\end{equation*}
(which is a bounded linear function of $p$ and $q$ in the given norms). With this choice of $\dot{c}$ the function $\dot{f}$ is
undetermined up to an additive constant; as usual, we take $\dot{f} = \ringf + e$ where $\ringf$ is the unique mean-value-zero
solution of the first two equations in \eqref{linear-system} and $e$ will be determined later. Notice that linear operator taking
$(p,q) \in L^2(\ENZ) \times H^{-1/2}(\partial D)$ to $\ringf \in S(\ENZ)$ is bounded.

We turn now to the third equation in \eqref{linear-system}. Remembering that the trace of $\dot{g}$ must match that of $\dot{f}$ at
$\partial D$, we see that it is to be solved with the Dirichlet boundary condition $\dot{g} = \dot{f}$ at $\partial D$. Since $\lambda_0$
is not a Dirichlet eigenvalue of $-\Delta$ in $D$, there is a unique solution; moreover it has the form
$$
\dot{g} = \ringg + e \psi_{d,\lambda_0}
$$
where $\ringg$ solves
\begin{equation} \label{ring-psi-pde}
\Delta \ringg + \la_0 \ringg + \psi_{d,\la_0} \dot{\mu} = r \quad \mbox{in $D$, with $\ringg = \ringf$ at $\partial D$.}
\end{equation}
With the benefit of foresight, we choose
\begin{equation} \label{mu-dot-choice}
\dot{\mu} = \frac{1}{|\ENZ| + \int_D \psi_{d,\lambda_0}^2} \,
\left[ \int_{\partial D} \ringf \partial_{\nu_D} \psi_{d,\lambda_0} + \int_D r (\psi_{d,\lambda_0} -1) \right]
\end{equation}
and
\begin{equation} \label{e-choice}
e = - \frac{1}{|\ENZ| + \int_D \psi_{d,\lambda_0}^2} \int_D \ringg \psi_{d,\lambda_0} \, .
\end{equation}
We note that $\dot{\mu}$ is a bounded linear functional of $(p,q,r) \in Y$, and it doesn't depend on $e$. Moreover, the
operator taking $(p,q,r) \in Y$ to $\ringg \in S(D)$ is linear and bounded, since $\ringg$ solves a Helmholtz-type PDE in $D$
whose source term $r - \dot{\mu} \psi_{d,\lambda_0}$ and Dirichlet data $\ringf$ are in $L^2(D)$ and $H^{1/2}(\partial D)$, each
depending linearly on $(p,q,r) \in Y$. Finally, since $\ringg$ depends linearly on $(p,q,r)$ and is independent of $e$, our
choice of $e$ is a bounded linear functional of $(p,q,r) \in Y$.

To complete the proof that our linear system is invertible, we must show that our choices \eqref{mu-dot-choice} and
\eqref{e-choice} assure the validity of the relations
\begin{equation} \label{two-crucial-relations}
\int_{\ENZ} \dot{f} + \int_D \dot{g} = 0 \quad \mbox{and} \quad \int_{\ENZ} \dot{f} + \int_D \dot{g}\psi_{d,\la_0} = 0 \, .
\end{equation}
To get the first, we multiply the PDE \eqref{ring-psi-pde} by $\psi_{d,\lambda_0} -1$, integrate over $D$, integrate by parts, and use the
Dirichlet boundary condition to get
$$
- \int_{\partial D} \ringf ( \partial_{\nu_D} \psi_{d,\lambda_0}) \, d \sh^1 - \lambda_0 \int_D \ringg +
\dot{\mu} \int_D (\psi_{d,\lambda_0}^2 - \psi_{d,\lambda_0}) \, dx = \int_D r(\psi_{d,\lambda_0} -1) \, dx \, .
$$
Since
$$
\int_D (\psi_{d,\lambda_0}^2 - \psi_{d,\lambda_0}) \, dx = |\ENZ| + \int_D \psi_{d,\lambda_0}^2 \, dx
$$
by the crucial consistency condition \eqref{consistency-first-version}, we see that \eqref{mu-dot-choice} is equivalent to
$$
\int_D \ringg \, dx = 0 \, .
$$
Since $\dot{f} = \ringf + e$ and $\dot{g} = \ringg + e \psi_{d,\lambda_0}$, we conclude that
$$
\int_{\ENZ} \dot{f} + \int_D \dot{g} = e \, (|\ENZ| +  \int_D \psi_{d,\lambda_0}) = 0,
$$
which gives the first equation in \eqref{two-crucial-relations}. As for the other, we have
$$
\int_{\ENZ} \dot{f} + \int_D \dot{g}\psi_{d,\lambda_0} = e \, |\ENZ| + \int_D (\ringg \psi_{d,\lambda_0} + e \psi_{d,\lambda_0}^2) \, dx \, ;
$$
evidently, our choice of $e$ in \eqref{e-choice} is exactly the one that makes this vanish.

We conclude, by the implicit function theorem, the existence of $\mu_\delta, f_\delta, g_\delta, c_\delta$ depending analytically
on $\delta$ in a (complex) neighborhood of $0$, such that $F(\delta, \mu_\delta, f_\delta, g_\delta, c_\delta) = 0$.
\medskip

\noindent {\sc Step 4} We now prove, using the specific structure of $F$, that in fact $c_\delta = 0$. Indeed, using
Green's theorem (but not the fact that $F=0$), we have
\begin{multline} \label{start-of-calculation}
\int_{\ENZ} \Bigl( \Delta f_\de + \la_0 + \de(\mu_\de + f_\de (\la_0 + \de \mu_\de)) + c_\de \Bigr) \\
 = - \int_{\partial D}\partial_{\nu_D} f_\de \,d\sh^1 + (\la_0 + \de \mu_\de) |\ENZ|  +
 \de( \la_0 + \de \mu_\de) \int_{\ENZ} f_\de + c_\de |\ENZ| \, .
\end{multline}
Adding and subtracting some terms, the right hand side can be rewritten as
\begin{multline*}
\int_{\partial D} \partial_{\nu_D}(\psi_{d,\la_0} + \de g_\de - f_\de) \,d \sh^1 -
 \int_{\partial D} \partial_{\nu_D}(\psi_{d,\la_0} + \de g_\de )\,d\sh^1  \\
 + (\la_0 + \de \mu_\de) |\ENZ|  + \de( \la_0 + \de \mu_\de) \int_{\ENZ} f_\de + c_\de |\ENZ| \, ,
\end{multline*}
which (applying Green's theorem) is the same as
\begin{multline*}
\int_{\partial D} \partial_{\nu_D}(\psi_{d,\la_0} + \de g_\de - f_\de) \,d \sh^1 - \int_D \Delta (\psi_{d,\la_0} + \de g_\de)\\
+ (\la_0 + \de \mu_\de) |\ENZ|  + \de( \la_0 + \de \mu_\de) \int_{\ENZ} f_\de + c_\de |\ENZ| \, .
\end{multline*}
Adding and subtracting some terms and using that $\Delta \psi_{d,\lambda_0} + \lambda_0 \psi_{d,\lambda_0} = 0$, the preceding
expression can be further rewritten as
\begin{multline*}
\int_{\partial D} \partial_{\nu_D}(\psi_{d,\la_0} + \de g_\de - f_\de) \,d \sh^1 -
 \delta \int_D (\Delta g_\de + (\la_0 + \de \mu_\de) g_\de + \mu_\de \psi_{d,\la_0} )\,dx +
 (\lambda_0 + \de \mu_\de) \int_D \psi_{d,\la_0}  \\
  + \de(\la_0 + \de\mu_\de) \int_D g_\de + (\la_0 + \de \mu_\de) |\ENZ|  + \de( \la_0 + \de \mu_\de) \int_{\ENZ} f_\de  + c_\de |\ENZ| \, .
\end{multline*}
Using now the fact that $F(\delta, \mu_\delta, f_\delta, g_\delta, c_\delta) = 0$, we conclude that
$$
(\lambda_0 + \de \mu_\de) \int_D \psi_{d,\la_0} + \de(\la_0 + \de\mu_\de) \int_D g_\de + (\la_0 + \de \mu_\de) |\ENZ|  +
\de( \la_0 + \de \mu_\de) \int_{\ENZ} f_\de  + c_\de |\ENZ| = 0 \, .
$$
Making use of the additional relations $\int_D g_\de + \int_{\ENZ} f_\de = 0$ and $|\ENZ| + \int_D \psi_{d,\lambda_0} =0$,
we finally conclude that $c_\delta |\ENZ| = 0$. Thus $c_\delta = 0$, as asserted.
\medskip

\noindent {\sc Step 5:} Remembering that $\mu_\delta$, $f_\delta$, and $g_\delta$ determine $u_\delta$ and $\lambda_\delta$ via
\eqref{e.proofansatz}, we have demonstrated the existence of an eigenpair $(u_\delta, \lambda_\delta)$ depending analytically on
$\delta$. The only remaining assertion of the theorem is that this is a simple eigenvalue, i.e. that the eigenspace of $\lambda_\de$ is
one-dimensional. This comes directly from the implicit function theorem, which tells us that
$z_\delta =(\mu_\delta, f_\delta, g_\delta, c_\delta)$ is the \emph{only} solution of $F(\delta, z) = 0$ near $z_0 = (\mu_0,f_0, g_0, 0)$
when $\delta$ is sufficiently small. If the eigenspace of $\lambda_\delta$ were multidimensional there would be
more than one solution of $F(\delta, z_\delta)$ near $(0,z_0)$; so the eigenspace is one-dimensional.
\end{proof}


\section{Accounting for dispersion and loss} \label{sec:accounting-for-dispersion-and-loss}

As noted in the Introduction, the dielectric permittivity of a material is generally a function of frequency
(this is known as \emph{dispersion}) and it is complex-valued (since waves decay as they propagate through materials).
Some key structural conditions are that
\begin{align}
& \eps(\omega) \ \mbox{is holomorphic in the upper half-plane,} \label{gen-prin-one}\\
& \eps(-\overline{\omega}) = \overline{\eps}(\omega) \ \mbox{for all $\omega \in \CC$, and} \label{gen-prin-two}\\
& \mbox{the imaginary part of $\eps(\omega)$ is nonnegative when $\omega$ is real and positive} \label{gen-prin-three}
\end{align}
(see e.g. Section 82 of \cite{landau-lifshitz}).
The class of all such functions is huge. When considering a particular material, however, a parsimonious framework is
needed, and for this purpose the \emph{Lorentz model} is often used. (In particular, \cite{liberal2016geometry} uses such a model to simulate silicon carbide as an ENZ material.) It has the form
\begin{equation} \label{lorentz-model-bis}
\eps(\omega,\gamma) = \eps_\infty \left( 1 + \frac{\omega_p^2}{\omega_0^2 - \omega^2 - i \omega \gamma} \right)
\end{equation}
where $\eps_\infty$, $\omega_p$, $\omega_0$, and $\gamma$ are nonnegative real numbers. (In discussing the dependence of this
function on $\omega$ with $\gamma$ held fixed, we shall sometimes omit the variable $\gamma$, writing $\eps(\omega)$
rather than $\eps(\omega, \gamma)$.) Viewed as a function of $\omega \in \CC$, this model has two poles in the
lower half-plane; to leading order as $\gamma \rightarrow 0$ they are at $-\frac{i}{2} \gamma \pm \omega_0$
(provided $\omega_0 \neq 0$). The Lorentz model is, roughly speaking, the simplest functional form
consistent with the general principles \eqref{gen-prin-one}--\eqref{gen-prin-three} (though it is
sometimes simplified further by taking $\omega_0 = 0$;
this is known as the Drude model).

\begin{figure}[h]
\begin{center}
\includegraphics[scale=.5]{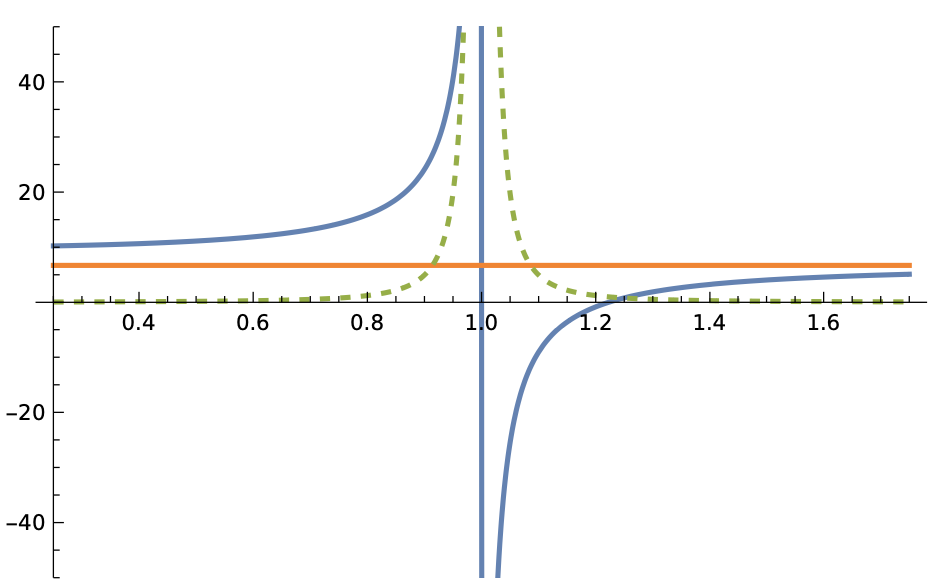}
\end{center}
\caption{The Lorentz model
$\eps(\omega,\gamma) = \eps_\infty \left( 1 + \frac{\omega_p^2}{\omega_0^2 - \omega^2 - i \omega \gamma} \right)$, graphed as
a function of $\omega/\omega_0$:
the solid blue curve is the real part of $\eps(\omega, \gamma)$, while the dotted green curve
is the imaginary part. The horizontal orange line shows the value of $\eps_\infty$. This figure was
produced using $\eps_\infty = 6.7$, $\omega_p/\omega_0 = 0.7$, and $\gamma/\omega_0 = .006$, consistent with
experimental data on silicon carbide near its resonance at frequency $\omega_0 = 2.38 \times 10^{13} \mbox{sec}^{-1}$
\cite{Spitzer}. This material system was used for the simulations in \cite{liberal2016geometry}.}
\label{fig:lorentz-model}
\end{figure}

Dispersion is more than just a fact of life -- it is in fact the \emph{reason} that ENZ materials exist. This
is especially easy to see for the Lorentz model. Indeed, in the lossless limit $\gamma = 0$ there is a
unique (real and positive) \emph{ENZ frequency}
\begin{equation} \label{enz-frequency}
\omega_* = \sqrt{\omega_p^2 + \omega_0^2}
\end{equation}
such that $\eps(\omega_*) = 0$. The presence of loss regularizes the singularity at $\omega = \omega_0$, but it leaves the
picture qualitatively intact: the real part of $\eps(\omega)$ vanishes at a $\gamma$-dependent real frequency near $\omega_*$.
The imaginary part of $\eps(\omega)$ is of course strictly positive when $\gamma>0$ and $\omega$ is real; however when
$\gamma$ is small it is mainly significant near $\omega_0$. (See Figure \ref{fig:lorentz-model}.)

The main result in this section, Theorem \ref{t.result-with-dispersion}, uses a Lorentz model for $\epsenz$
(though as we discuss in Remark \ref{rmk:more-general-dispersion} our method applies more generally).
We do not use a specific model for $\eps_D$; rather, we assume only that
\begin{align}
& \mbox{$\eps_D(\omega)$ is real-valued when $\omega$ is real;} \label{epsD-real}\\
& \mbox{$\eps_D(\omega_*)$ is positive, and $\eps_D$ is analytic in a neighborhood of $\omega_*$; and} \label{epsD-positive} \\
& \mbox{for real-valued $\omega$ near $\omega_*$, $\frac{d}{d \omega} [\omega^2 \eps_D(\omega)] > 0$}. \label{epsD-monotone}
\end{align}
The first condition says that the material in region $D$ has negligible loss at frequencies near $\omega_*$. (This
was assumed in \cite{liberal2016geometry}.) The second is very routine. The third condition is actually
satisfied by any physical material,
since when loss is negligible it is known that $\frac{d}{d \omega} [\omega \eps(\omega)] > 0$ when $\omega$ is
real and positive (see e.g. Section 80 of \cite{landau-lifshitz}).

\begin{theorem} \label{t.result-with-dispersion}
Let $\epsenz = \epsenz(\omega, \gamma)$ have the form \eqref{lorentz-model} for some $\omega_0 \geq 0$ and $\omega_p > 0$
(which will be held fixed), and let $\omega_*$ be the associated ENZ frequency \eqref{enz-frequency}. Suppose further that
$\eps_D$ satisfies \eqref{epsD-real}--\eqref{epsD-monotone}, that
\begin{equation} \label{lambda*-defn}
\lambda_* : = \frac{1}{c^2} \omega_*^2 \eps_D(\omega_*)
\end{equation}
is not a Dirichlet eigenvalue of $-\Delta$ in $D$, and that $\lambda_*$ satisfies the crucial consistency condition
\begin{equation} \label{consistency-condition-met}
|\ENZ| + \int_D \psi_{d,\lambda_*} \, dx = 0
\end{equation}
(which is \eqref{consistency-first-version} with $\lambda_0$ replaced by $\lambda_*$). Then there is an analytic function $\omega(\gamma)$
defined in a neighborhood of $0$ such that $\omega(0)= \omega_*$ and
\begin{equation} \label{eqn-to-be-solved}
\lambda_{\epsenz (\omega(\gamma),\gamma)/\eps_D(\omega(\gamma))} = \omega^2(\gamma) c^{-2} \eps_D(\omega(\gamma)) ,
\end{equation}
where $\lambda_\delta$ is the function supplied by Theorem \ref{t.main1} with $\lambda_0$ replaced by $\lambda_*$. It follows
that \eqref{helmholtz-D}--\eqref{helmholtz-ENZ} has a one-dimensional solution space when $\omega = \omega(\gamma)$, spanned by the function $u_\delta$
provided by Theorem \ref{t.main1} with $\delta = \epsenz (\omega(\gamma),\gamma)/\eps_D(\omega(\gamma))$. The
value of $\omega'(0)$ can be expressed in terms of
\begin{align}
a_1 := & \ \partial_\omega \epsenz(\omega_*,0) \label{defn-a1}\\
a_2 := & \ \frac{1}{i} \partial_\gamma \epsenz(\omega_*,0) \label{defn-a2}\\
a_3 : = & \ \partial_\omega (\omega^2 \eps_D (\omega)) |_{\omega = \omega_*} \label{defn-a3}
\end{align}
(all of which are easily seen to be positive real numbers) by
\begin{equation} \label{omega-prime-formula}
\omega'(0)  = -i \, \frac{a_2}{a_1 + a_3 c^{-2} \eps_D(\omega_*) |\lambda'(0)|^{-1}},
\end{equation}
where $\lambda'(0)$ is given by \eqref{e.lambda1def}.
\end{theorem}

Before giving the proof, let us discuss a key consequence of this result. When designing a resonator, it is natural to use
materials with relatively little loss, so the value $\gamma$ should be small and
$\omega(\gamma) \approx \omega(0) + \omega'(0) \gamma = \omega_* + \omega'(0) \gamma$.
Since $\omega'(0)$ is purely imaginary and $\gamma$ is positive,
we see that the real part of $\omega(\gamma)$ (which is, physically speaking, the resonant frequency) is very near the
ENZ frequency $\omega_*$ (the difference is at most of order $\gamma^2$). We also see that the imaginary part of
$\omega(\gamma)$ (which controls the quality factor of the resonance -- in other words the rate at which it decays)
depends on the shape of $\Omega$ only through $|\lambda'(0)|$, and that the quality factor is optimized
(the decay rate is minimized) by choosing the shape of $\Omega$ so that $|\lambda'(0)|$ is as small as possible.

\begin{proof}[Proof of Theorem \ref{t.result-with-dispersion}] By the implicit function theorem, it suffices to show
that when we calculate $\omega'(0)$ formally by differentiating \eqref{eqn-to-be-solved}, the calculation succeeds (without
dividing by $0$). Remembering that $\epsenz(\omega_*,0)=0$, differentiation with respect to $\gamma$ at $\gamma = 0$ gives
$$
\lambda'(0) \left[ \frac{a_1 \omega'(0)}{\eps_D(\omega_*)} + \frac{ia_2}{\eps_D(\omega_*)}\right] =  a_3 c^{-2} \omega'(0).
$$
Solving for $\omega'(0)$ gives
$$
\omega'(0) = -i \frac{a_2}{a_1 - a_3 c^{-2} \eps_D(\omega_*) (\lambda'(0))^{-1}}.
$$
Since we know from \eqref{e.lambda1def} that $\lambda'(0)$ is a negative real number, the preceding expression is
equivalent to \eqref{omega-prime-formula}.
\end{proof}

\begin{remark} \label{rmk:more-general-dispersion}
While we have assumed, for simplicity, that $\epsenz$ is given by a Lorentz model, our method is clearly
also applicable in other settings. Its key requirements are that (i) $\epsenz = \epsenz(\omega,\gamma)$ be a function of
the frequency $\omega$ and a single (scalar) loss parameter $\gamma$, and that (ii) its partial derivatives at $\gamma=0$,
$\omega = \omega_*$ be such that $a_1$ and $a_2$ are positive real numbers. Suppose, for example, that the permittivity of
the ENZ material has the form
$$
\eps(\omega) = \eps_\infty \left( 1 + \sum_{j=1}^N \frac{(\omega_p^j)^2}{(\omega_0^j)^2 - \omega^2 - i \omega \gamma^j} \right)
$$
for some (positive, real) constants $\omega_p^j$, $\omega_0^j$, and $\gamma^j$, ordered so that
$\omega_0^1 < \ldots < \omega_0^N$. By the discussion associated with Figure \ref{fig:lorentz-model}, such a material
has an ENZ frequency $\omega_*^j$ (defined as a root of $\eps(\omega)=0$ when $\gamma^1, \ldots, \gamma^N$
are all set to $0$) between $\omega_0^j$ and $\omega_0^{j+1}$ for each $j=1,\ldots,N-1$. To get a resonance near
$\omega_*^1$ (for example), it is natural to use
$$
\epsenz(\omega,\gamma) =
\eps_\infty \left( 1 + \sum_{j=1}^N \frac{(\omega_p^j)^2}{(\omega_0^j)^2 - \omega^2 - i \omega \gamma \hat{\gamma}^j} \right)
$$
with $\hat{\gamma}^j = \gamma^j/\gamma^1$. One easily checks that $a_1$ and $a_2$ are positive, so
our implicit-function-theorem-based argument is applicable. (However our result is local: it gives a resonant frequency
$\omega(\gamma)$ for $\gamma$ near $0$. The argument does not show that $\omega(\gamma)$ is defined even for
$\gamma = \gamma^1.$)
\end{remark}

\section{The optimal design problem} \label{sec:the-optimal-design-problem}

Theorem \ref{t.result-with-dispersion} proves the existence of a resonance at (complex) frequency $\omega(\gamma)$ when the loss
parameter $\gamma$ is sufficiently close to $0$. The theorem's hypotheses involve the area of $\Omega$, but they are otherwise
independent of its shape. However, according to eqn. \eqref{omega-prime-formula} the quality of the resonance \emph{does} depend on
the shape of $\Omega$. Therefore it is natural to ask how $\Omega$ should be chosen so as to \emph{optimize} the resonance.
Theorem \ref{t.result-with-dispersion} shows that, to leading order in $\gamma$, this amounts to asking what shape
minimizes $|\lambda'(0)|$.

The function $\lambda(\delta)$ was introduced in Section \ref{sec:analysis-without-dispersion}, where our notation was
$\lambda(0) = \lambda_0$ and $\lambda'(0) = \lambda_1$. The analysis in Section \ref{sec:accounting-for-dispersion-and-loss}
used a particular choice of $\lambda_0$, which we called $\lambda_*$. However our optimal design problem can be considered
for any choice of $\lambda_0$. Therefore we revert in this section to the notation of Section \ref{sec:analysis-without-dispersion}.

In considering this optimal design problem, we will be holding $D$ and $\lambda_0$ fixed. It follows from the consistency condition
that $|\ENZ|$ is also fixed. Recalling from \eqref{e.lambda1def} that
$$
\lambda_1 = - \frac{1}{|\ENZ| + \int_D \psi_{d,\la_0}^2 \, dx} \int_D |\nabla \phi_1|^2 \, dx
$$
and observing that the expression in front of the integral is being held fixed, we see that the goal of our
optimal design problem is to minimize the value of
\begin{equation} \label{to-be-optimized}
\frac12 \int_{\ENZ}|\nabla \phi_1|^2 \, dx \, ,
\end{equation}
where $\phi_1$ solves \eqref{e.neumann}, which we repeat for the reader's convenience here:
\begin{equation} \label{e.phi1def.def}
\begin{aligned}
-\Delta \phi_1 &= \lambda_0 \quad \mbox{ in } \ENZ\\
\partial_{\nu_\Omega} \phi_1 &= 0 \quad \mbox{ on } \partial \Omega\\
\partial_{\nu_D} \phi_1 &= \partial_{\nu_D}\psi_{d,\la_0} \quad \mbox{ on } \partial D \, .
\end{aligned}
\end{equation}
(Since this is a pure Neumann problem, the data must be consistent; this is assured by the consistency
condition \eqref{consistency-first-version}, as we showed in Section \ref{subsec:leading-order-terms}.
The solution is only unique up to a constant, but the value of \eqref{to-be-optimized} is independent of this constant.)

It is a standard fact that \eqref{to-be-optimized} has a variational characterization:
\begin{equation} \label{variational-characterization}
-\frac12 \int_{\ENZ}|\nabla \phi_1|^2 \, dx = \min_{w \in H^1(\ENZ)}
\int_{\ENZ} \frac12 |\nabla w|^2 - \lambda_0 w  \, dx +
\int_{\partial D} (\partial_{\nu_D} \psi_{d,\lambda_0} )w \, d\sh^1 \, ,
\end{equation}
and that $\phi_1$ is optimal for RHS of \eqref{variational-characterization}.
(To explain the sign of the boundary term in the variational principle, we note that $\nu_D$ is the \emph{inward}
unit normal to $\partial (\ENZ)$ at $\partial D$.)
Our optimal design problem can thus be restated as
\begin{equation} \label{max-min}
\sup_{\Omega \supset \overline{D}} \, \inf_{w \in H^1(\ENZ)}
\int_{\ENZ}  \frac12 |\nabla w|^2 - \lambda_0 w  \, dx +
\int_{\partial D} (\partial_{\nu_D} \psi_{d,\lambda_0} )w \, d \sh^1 \, ,
\end{equation}
where it is understood that $\Omega$ ranges over Lipschitz domains. It might seem that the optimization over $\Omega$
should be subject to a constraint on $|\ENZ|$, in view of the consistency condition \eqref{consistency-first-version}.
Actually no such constraint is needed, since if the consistency condition is violated then
the minimization over $w$ takes the value $-\infty$. (In \eqref{max-min} and throughout this section, we write
$\inf$ and $\sup$ rather than $\min$ and $\max$ when we do not mean to claim that the optimum is achieved.)

We have two main results on this optimal design problem:
\begin{itemize}
\item In Section \ref{subsec:round-D} we show that if $D$ is a ball then the optimal $\Omega$ is a concentric ball.
\item In Section \ref{subsec:convex-relaxation} we study a convex relaxation of \eqref{max-min}, which is certainly an upper
bound but which is conjecturally equivalent to the unrelaxed problem.
\end{itemize}
The relationship between our relaxation of \eqref{max-min} and the unrelaxed problem is discussed in
Section \ref{subsec:relaxation}. As we explain there, our relaxation has a physical interpretation involving
homogenization. This use of homogenization is similar to the introduction of composite materials in
compliance optimization problems with design-independent loading, as studied for example in \cite{Allaire}.
While this interpretation of our relaxation has yet to be justified rigorously for \eqref{max-min}, it has
been fully justified for compliance optimization problems with design-independent loading.


\subsection{Optimality of a ball for round $D$} \label{subsec:round-D}

\begin{theorem} \label{t.round-D}
If $D$ is a ball, then $|\lambda_1|$ is minimized by taking $\Omega$ to be a concentric ball. (Its radius is
determined by $D$ and $\lambda_0$ through the consistency condition.) Moreover, this optimum is unique: no other $\Omega$
can do as well.
\end{theorem}
\begin{proof}
It suffices to consider the case when $D$ is the unit disk, since the general case is easily reduced to this one by
translation and scaling. The function $\psi_{d,\lambda_0}$ is then radial and quite explicit:
\begin{equation} \label{formula-for-psi-d}
\psi_{d,\lambda_0}(r) = \frac{J_0(\lambda_0 r)}{J_0(\lambda_0)}, \quad r \in (0,1) \, ,
\end{equation}
where as usual $J_0$ is the zeroth order cylindrical Bessel function of the first kind. Since $\lambda_0$ is
not an eigenvalue of the Laplacian in the unit disc, this is well-defined ($J_0(\lambda_0) \neq 0$).

Let $A_0 = |\ENZ|$ be the area of the ENZ shell. Its value is available from the consistency condition for the existence
of $\phi_1$:
$$
\lambda_0  A_0 = \int_{\partial D}\partial_{\nu_D} \psi_{d,\lambda_0}\,d\sh^1 \, ,
$$
which in view of \eqref{formula-for-psi-d} gives
$$
A_0 = 2\pi \frac{J_0^{\prime}(\lambda_0)}{J_0(\lambda_0)} \, .
$$
(As noted in earlier sections, we need $\int_{\partial D}\partial_{\nu_D} \psi_{d,\lambda_0}\,d\sh^1 > 0$
in order that the area of $\ENZ$ be positive. This is a condition on $\lambda_0$, which reduces in the present
setting to $J_0^{\prime}(\lambda_0)/J_0(\lambda_0) > 0$.)

Our claim is that the optimal $\Omega$ is a ball centered at the origin with area $|D| + |\ENZ| = \pi + A_0$.
Let us call this domain $\Omega_0$; it is the ball whose radius $r_0$ satisfies $\pi r_0^2 = \pi + A_0$. The function
$\phi_1$ associated with $\Omega_0$ is easily made explicit. Since it is clearly radial, we may write
$\phi_1 = \phi_1(r)$, so that the boundary value problem \eqref{e.phi1def.def} becomes
\begin{align*}
- \phi_1^{\prime\prime}(r) - \frac{\phi_1^\prime(r)}{r} &= \lambda_0 \quad \mbox{ for } r \in (1,r_0)\\
\phi_1^\prime(r_0) &= 0 \\
\phi_1^\prime (1) &= \psi_{d,\lambda_0}^\prime(1) = \lambda_0 \frac{J_0^{\prime}(\lambda_0)}{J_0(\lambda_0)} \,.
\end{align*}
The solution is unique up to an additive constant. The general solution of the ODE is
$\phi_1(r) = b + c \log r - \lambda_0 \frac{r^2}{4}$, and the boundary condition at $r=r_0$ gives $c=\lambda_0 r_0^2/2$.
(The boundary condition at $r=1$ gives no additional information; it is automatically satisfied, as a consequence of
the consistency condition.)

While the constant $b$ is arbitrary, it is convenient to choose it so that $\phi_1(r_0) = 0$. The resulting (now fully
determined) function has the property that $\phi_1(r) < 0$ for $1 \leq r < r_0$. (Indeed, $\phi_1$ is strictly concave and
$\phi_1'(r_0) = 0$, so it is an increasing function on this interval and it vanishes at $r=r_0$.) This implies in particular
that
\begin{equation} \label{e.phi1sign}
\frac12 |\nabla \phi_1|^2 - \lambda_0 \phi_1 =
\frac12 (\phi_1^\prime)^2 - \lambda_0 \phi_1 > 0 \quad \mbox{ for } r \in (1,r_0)\,.
\end{equation}

To demonstrate the optimality of $\Omega_0$, we shall use the extension of $\phi_1$ by $0$,
$$
\tilde{\phi}_1 = \left\{
\begin{array}{cl}
\phi_1 (r) & \mbox{for $1 \leq r \leq r_0$}\\
0 & \mbox{for $r \geq r_0$,}
\end{array}
\right.
$$
as a test function in the variational principle that characterizes $\lambda_1$.

Let $\hat{\Omega}$ be a competitor to $\Omega_0$; in other words, let $\hat{\Omega} \subset \RR^2$ be a bounded, open set
with locally Lipschitz boundary that contains $\overline{D}$ and satisifes $|\hat{\Omega} \setminus \overline{D}|= A_0,$.
If $\hat{\phi}_1$ is the solution of \eqref{e.phi1def.def} with $\hat{\Omega}$ in place of $\Omega$, then the
variational principle \eqref{variational-characterization} gives
\begin{align}
- \frac12 \int_{\hat{\Omega} \setminus \overline{D}} |\nabla \hat\phi|^2\,dx &=
 \int_{\hat{\Omega} \setminus \overline{D}} \frac12 |\nabla \hat{\phi}|^2 - \lambda_0 \hat{\phi} \,dx +
 \int_{\partial D} (\partial_{\nu_D}\psi_{d,\lambda_0})\hat{\phi}\,d\sh^1 \nonumber\\
&= \min_{w \in H^1(\hat{\Omega} \setminus \overline{D})} \int_{\hat{\Omega} \setminus \overline{D}}
\frac12 |\nabla w|^2 - \lambda_0 w \,dx + \int_{\partial D} (\partial_{\nu_D}\psi_{d,\lambda_0})w\,d\sh^1 \nonumber\\
&\leq \int_{\hat{\Omega} \setminus \overline{D}} \frac12 |\nabla \tilde{\phi}_1|^2 - \lambda_0 \tilde{\phi}_1  \,dx +
\int_{\partial D} (\partial_{\nu_D}\psi_{d,\lambda_0}) \tilde{\phi}_1\,d\sh^1 \, . \label{intermediate-step}
\end{align}
Since $\frac12 |\nabla \tilde{\phi}_1|^2 - \lambda_0 \tilde{\phi}_1$ vanishes outside $\Omega_0$ and is positive in
$\Omega_0 \setminus \overline{D}$,
\begin{align} \label{key-to-uniqueness}
\int_{\hat{\Omega} \setminus \overline{D}} \frac12 |\nabla \tilde{\phi}_1|^2 - \lambda_0 \tilde{\phi}_1  \,dx &=
\int_{(\hat{\Omega} \setminus \overline{D}) \cap \Omega_0} \frac12 |\nabla \phi_1|^2 - \lambda_0 \phi_1  \, dx \nonumber \\
&\leq  \int_{\Omega_0 \setminus \overline{D}} \frac12 |\nabla \phi_1|^2 - \lambda_0 \phi_1  \,dx \, .
\end{align}
Combining this with \eqref{intermediate-step} gives
\begin{align*}
- \frac12 \int_{\hat{\Omega} \setminus \overline{D}} |\nabla \hat\phi|^2\,dx &\leq
\int_{\Omega_0 \setminus \overline{D}} \frac12 |\nabla \phi_1|^2 - \lambda_0 \phi_1  \,dx +
\int_{\partial D} (\partial_{\nu_D}\psi_{d,\lambda_0})\phi_1\,d\sh^1 \\
&= - \frac12 \int_{\Omega_0 \setminus \overline{D}} |\nabla \phi_1|^2\,dx \,
\end{align*}
where in the final step we used \eqref{variational-characterization}. This confirms the optimality of $\Omega_0$. To see
its uniqueness, we recall that $\frac{1}{2} |\nabla \phi_1|^2 - \lambda_0 \phi_1$ is strictly positive in
$\Omega_0 \setminus \overline{D}$.
Therefore equality holds in \eqref{key-to-uniqueness} only when $\hat{\Omega} \setminus \overline{D}$
includes the entire domain $\Omega_0 \setminus \overline{D}$. Since both sets have area $A_0$, it follows that
$\hat{\Omega} \setminus \overline{D}= \Omega_0 \setminus \overline{D}$, whence $\hat{\Omega} = \Omega_0$.
\end{proof}

\begin{remark}
The preceding argument is simple, but perhaps a bit mysterious. The next section offers a convex-optimization-based
perspective on our optimal design problem. In general, for a convex variational problem, if one can guess the optimal
test function, then there is usually a simple proof that the guess is right, obtained by using a solution of the
dual problem. As we shall show in Proposition \ref{classical-solutions-of-relaxed-problem},
this is indeed the character of the argument just presented.
\end{remark}

\subsection{A convex relaxation} \label{subsec:convex-relaxation}

We turn now to the max-min problem \eqref{max-min}, when $D$ is any simply-connected Lipschitz domain.
We start by making some minor adjustments:
\begin{itemize}
\item As noted at the beginning of Section \ref{subsec:preliminaries-etc}, we do not want to assume that $\Omega$ is simply
connected. However we want $\Omega$ to be a \emph{bounded} domain, and it is therefore natural to introduce the restriction
that $\Omega$ be a subset of some fixed region $B$ that contains $\overline{D}$.

\item The function $\partial_{\nu_D} \psi_{d,\lambda_0}$ appears in the final term of \eqref{max-min}, because the PDE
for $\phi_1$ is driven by this source term at $\partial D$. However, the analysis in this section applies equally when
this function is replaced by any $f \in H^{-1/2}(\partial D)$ such that
\begin{equation} \label{condition-on-f}
\int_{\partial D} f \, d \sh^1 > 0 .
\end{equation}
To emphasize this, throughout the present section our source term will be $f$ rather than $\partial_{\nu_D} \psi_{d,\lambda_0}$.

\item When we replace $\partial_{\nu_D} \psi_{d,\lambda_0}$ by $f$ in the PDE \eqref{e.phi1def.def} defining $\phi_1$, the condition
for existence of a solution becomes
$$
\lambda_0 |\ENZ| = \int_{\partial D} f\, d \sh^1 .
$$
Obviously $B$ must be large enough to contain $\Omega$, so we require that
\begin{equation} \label{B-large-enough}
|B \setminus \overline{D}| > \frac{1}{\lambda_0} \int_{\partial D} f\, d \sh^1 .
\end{equation}
\end{itemize}

\noindent Taking these adjustments into account, our goal is to understand
\begin{equation} \label{max-min-B}
m := \sup_{\Omega \ {\rm s.t.}\, \overline{D} \subset \Omega \subset B} \ \inf_{w \in H^1(\ENZ)}
\int_{\ENZ} \frac12 |\nabla w|^2 - \lambda_0 w \, dx +
\int_{\partial D} f w \, d \sh^1 \, ,
\end{equation}
with the unspoken convention that $\Omega$ ranges over Lipschitz domains.
It is convenient to write this differently, in terms of the \emph{characteristic function} of $\Omega$,
viewed as a function on $B \setminus \overline{D}$ that takes only the values $0$ and $1$ (outside and inside $\Omega$
respectively):
\begin{equation} \label{max-min-bis}
m = \sup_{
\substack{\chi(x) \in \{0,1\} \, {\rm for} \, x \in B \setminus \overline{D}\\
\chi = 1 \, {\rm at} \, \partial D}
} \ \inf_{w \in H^1(B \setminus \overline{D})}
\int_{B \setminus \overline{D}} \chi(x) \left( \frac12 |\nabla w|^2 - \lambda_0 w \right) \, dx +
\int_{\partial D} f w \, d \sh^1 \, .
\end{equation}
Our convex relaxation of the optimal design problem is obtained by replacing the characteristic
function $\chi$ (which takes only the values $0$ and $1$) by a density $\theta$ (which takes any value
$0 \leq \theta \leq 1$). Since enlarging the class of test
functions in a maximization can only increase the value of the maximum, it is obvious that
\begin{equation} \label{max-min-relaxed}
m \leq m_{\rm rel} = \sup_{0 \leq \theta(x) \leq 1} \ \inf_{w \in H^1(B \setminus \overline{D})}
\int_{B \setminus \overline{D}} \theta(x) \left( \frac12 |\nabla w|^2 - \lambda_0 w \right) \, dx +
\int_{\partial D} f w \, d \sh^1 \, .
\end{equation}
There is reason to think that $m = m_{\rm rel}$, as we shall explain in Section \ref{subsec:relaxation}. For now,
however, we focus on the relaxed problem \eqref{max-min-relaxed}.

It might seem strange that in formulating the relaxed problem we have kept no remnant of the condition that
$\chi = 1$ at $\partial D$. The reason is that if $\chi = 0$ near a part of $\partial D$ where
$ f \neq 0$, then the min over $w$ is $-\infty$ (by considering test functions $w$
supported in the region where $\chi = 0$). So we believe that the value of \eqref{max-min-bis} is not
changed by dropping the constraint that $\chi = 1$ at $\partial D$.

Equation \eqref{max-min-relaxed} defines $m_{\rm rel}$ as the sup-inf of
\begin{equation} \label{definition-of-L}
L(w,\theta) = \int_{B \setminus \overline{D}} \theta(x) \left( \frac12 |\nabla w|^2 - \lambda_0 w \right) \, dx +
\int_{\partial D} f w \, d \sh^1
\end{equation}
We note that $L$ is linear in $\theta$ and convex in $w$, so the inf over $w$ in \eqref{max-min-relaxed} is a concave
function of $\theta$ and the sup-inf can be viewed as maximizing a concave function of $\theta$ subject to the convex
constraint $ 0 \leq \theta(x) \leq 1$. Formally, at least, the associated \emph{dual problem} is obtained by replacing
sup-inf by inf-sup:
\begin{equation} \label{dual-as-min-max}
\inf_{w \in H^1(B \setminus \overline{D})} \ \sup_{0 \leq \theta(x) \leq 1}
\int_{B \setminus \overline{D}} \theta(x) \left( \frac12 |\nabla w|^2 - \lambda_0 w \right) \, dx +
\int_{\partial D} f w \, d \sh^1 \, .
\end{equation}
It is obvious that
$$
\sup_{0 \leq \theta(x) \leq 1}
\int_{B \setminus \overline{D}} \theta(x) \left( \frac12 |\nabla w|^2 - \lambda_0 w \right) \, dx =
\int_{B \setminus \overline{D}} \left( \frac12 |\nabla w|^2 - \lambda_0 w \right)_+ \, dx
$$
with the notation $z_+ = \max\{z,0\}$, so the formal dual is equivalent to
\begin{equation} \label{simplified-dual}
\inf_{w \in H^1(B \setminus \overline{D})}
\int_{B \setminus \overline{D}} \left( \frac12 |\nabla w|^2 - \lambda_0 w \right)_+ \, dx +
\int_{\partial D} f w \, d \sh^1 \, .
\end{equation}
(We will show in due course that this infimum is achieved; but we note here that the functional tends to
infinity when $w=c$ is constant and $c \rightarrow \pm \infty$, as an easy consequence of \eqref{condition-on-f} and
\eqref{B-large-enough}.)

The following theorem justifies the preceding formal calculation; in particular, it shows that the optimal values of our primal
and dual problems are the same, and it proves the existence of an optimal $\theta$ for
\eqref{max-min-relaxed} and an optimal $w$ for \eqref{simplified-dual}.

\begin{theorem} \label{t.relaxed-design-thm}
Let
$$
{\mathcal B} = \{ \theta \in L^2 (B \setminus \overline{D}) \mbox{ such that } 0 \leq \theta(x) \leq 1 \mbox{ a.e.} \} ,
$$
and observe that $L(w,\theta)$ (defined by \eqref{definition-of-L}) is well-defined and finite for
$\theta \in {\mathcal B}$ and $w \in H^1(B \setminus \overline{D})$. Then

\begin{enumerate}
\item[(a)] there is a \emph{saddle point} $\overline{w} \in H^1(B \setminus \overline{D})$ and
$\overline{\theta} \in {\mathcal B}$, in other words a pair such that
$$
L(\overline{w},\theta) \leq L(\overline{w}, \overline{\theta}) \leq L(w,\overline{\theta})
$$
for all $\theta \in {\mathcal B}$ and $w \in H^1(B \setminus \overline{D})$; moreover
\item[(b)] the sup-inf \eqref{max-min-relaxed} and the inf-sup \eqref{dual-as-min-max} have the same value,
namely $L(\overline{w}, \overline{\theta})$.
\end{enumerate}
\end{theorem}

\begin{proof}
We will apply Proposition 2.4 from Chapter 6 of \cite{ET}. The overall framework of that chapter
involves a functional $L(w,\theta)$ which is defined (and finite) as $w$ and $\theta$ range over closed convex subsets
of reflexive Banach spaces. This framework applies to our example, since $w$ ranges over the entire space
$H^1(B \setminus \overline{D})$ and $\theta$ ranges over ${\mathcal B}$, which is a closed convex subset of
$L^2(B \setminus \overline{D})$. Chapter 6 of \cite{ET} needs the additional structural conditions that
$\theta \rightarrow L(w,\theta)$ be concave and upper semicontinuous as a function of $\theta$ when $w$
is held fixed, and that $w \rightarrow L(w,\theta)$ be convex and lower semicontinuous as a function of
$w$ when $\theta \in {\mathcal B}$ is held fixed. Our example meets these requirements.

Proposition 2.4 of \cite{ET} has two further hypotheses, namely that
\begin{enumerate}
\item[(i)] the constraint set ${\mathcal B}$ is bounded, and

\item[(ii)] there exists $\theta_0 \in {\mathcal B}$ such that
\begin{equation} \label{coercivity-of-L}
\lim_{\|w\| \rightarrow \infty} L(w,\theta_0) = \infty \ .
\end{equation}
\end{enumerate}
While (i) is valid in our situation, (ii) is not, since it fails when we restrict attention to constant $w$. To deal
with this difficulty, we will proceed in two steps. In the first we restrict $\theta$ to lie in the smaller constraint set
\begin{equation} \label{B-tilde-defn}
\tilde{\mathcal B} =  {\mathcal B} \cap \left\{ \theta \mbox{ such that } \lambda_0 \int_{B \setminus \overline{D}} \theta \, dx = \int_{\partial D} f \, d \sh^1  \right\} \, ,
\end{equation}
which is nonempty by \eqref{B-large-enough}. For such $\theta$, $L(w,\theta)$
has the property that $L(w,\theta) = L(w+c, \theta)$ for any constant $c$; therefore it
can be viewed as being defined for all $w \in H^1/\R$. (Here and in the rest of this proof, we use
$H^1/\R$ as shorthand for the space $H^1 (B \setminus \overline{D})/\R$.) In Step $1$ we will show
that the saddle point result from \cite{ET} applies when $\theta$ ranges over $\tilde{\mathcal B}$ and
$w$ ranges over $H^1/\R$. Then in Step $2$ we will use this result to prove the theorem.
\medskip

\noindent {\sc Step 1.} To apply the proposition from \cite{ET}, it suffices to show that
\eqref{coercivity-of-L} is valid when we choose $\theta_0 \in \tilde{\mathcal B}$ to have a
positive lower bound (for example, we could choose it to be constant), if we view
$w \rightarrow L(w,\theta_0)$ as a function on $H^1/\R$. This is standard; indeed, we may take each equivalence class
in $H^1/\R$ to be represented by a function with $\int_{B \setminus \overline{D}} w \, dx = 0$.
By Poincar\'{e}'s inequality, we may take the norm on $H^1/\R$ to be
$\| \nabla w \|_{L^2(B \setminus \overline{D})}$. By the trace theorem and
Poincar\'{e}'s inequality, the terms in $L$ that are linear in $w$ are bounded by a constant times $\|w\|$, whereas
$$
\int_{B \setminus \overline{D}} \theta_0 |\nabla w|^2 \, dx \geq c \|w\|^2
$$
where $c$ is a lower bound for $\theta_0$. The quadratic term dominates when $\|w\|$ is large enough, so \eqref{coercivity-of-L} holds.
With the notation
$$
A_0 = \frac{1}{\lambda_0} \int_{\partial D} f \, d \sh^1 \,
$$
we conclude (applying the result from \cite{ET}) that
\begin{equation} \label{step-one-conclusion}
\sup_{
\substack{0 \leq \theta \leq 1\\
\int \theta(x) \, dx = A_0}
} \inf_{w \in H^1/\R} L(w,\theta) =
\inf_{w \in H^1/\R} \sup_{
\substack{0 \leq \theta \leq 1\\
\int \theta(x) \, dx = A_0}
} L(w,\theta) ;
\end{equation}
that there exist $\tilde{\theta}$ (satisfying $0 \leq \tilde{\theta} \leq 1$ and $\int \tilde{\theta} \, dx = A_0$)
and $\tilde{w}$ (in $H^1/\R$) satisfying
$$
L(\tilde{w},\theta) \leq L(\tilde{w}, \tilde{\theta}) \leq L(w,\tilde{\theta})
$$
for all $w \in H^1/\R$ and all $0\leq \theta \leq 1$ satisfying
$\int_{B \setminus \overline{D}} \theta \, dx = A_0$;
and that the value of \eqref{step-one-conclusion} is $L(\tilde{w},\tilde{\theta})$.
\bigskip

\noindent {\sc Step 2.} The desired saddle point $(\overline{w}, \overline{\theta})$ will be $(w_1,\tilde{\theta})$, where
$w_1$ is a well-chosen representative of $\tilde{w}$. To get started, let us examine the relationship between
$\tilde{w}$ and $\tilde{\theta}$, using the fact that $\tilde{\theta}$ maximizes $L(\tilde{w},\theta)$ over all
$\theta \in \tilde{B}$. Since $L(w,\theta) = L(w+c, \theta)$ when $c$ is constant and $\theta \in \tilde{B}$, we may work
with any representative $w_0 \in H^1(B \setminus \overline{D})$ of $\tilde{w}$. Evidently, $\tilde{\theta}$ achieves
\begin{equation} \label{optimality-of-theta-tilde}
\sup_{
\substack{0 \leq \theta \leq 1\\
\int \theta(x) \, dx = A_0}
}
\int_{B \setminus \overline{D}} \theta(x) \left( \frac12 |\nabla w_0|^2 - \lambda_0 w_0 \right) \, dx +
\int_{\partial D} f w_0 \, d\sh^1 \, .
\end{equation}
Since $\theta$ doesn't enter the boundary term, we shall be focusing in what follows on the bulk term.
To understand what conclusions we can draw from the optimality of $\tilde{\theta}$, let us assume for a moment that
$\frac12 |\nabla w_0|^2 - \lambda_0 w_0$ has no level sets with positive measure. Then there is a
unique $z_0 \in \R$ such that
$$
\left| \left\{ x \in B \setminus \overline{D} \, : \, \frac12 |\nabla w_0|^2 - \lambda_0 w_0 > z_0 \right\} \right| = A_0
$$
and $\tilde{\theta}$ must be the characteristic function of this set. In general, however, we must
allow for the possibility that $\frac12 |\nabla w_0|^2 - \lambda_0 w_0 $ has level sets with positive
measure. To deal with this, let
$$
g(z) = \left| \left\{ x \in B \setminus \overline{D} \, : \, \frac12 |\nabla w_0|^2 - \lambda_0 w_0 > z \right\} \right| \, ,
$$
which is a monotone (but possibly discontinuous) function of $z \in \R$. We can then consider two cases:
\begin{enumerate}
\item[(i)] If there exists $z_0$ such that $g(z_0) = A_0$, then $\tilde{\theta}$ must be the characteristic function
of the set where $\frac12 |\nabla w_0|^2 - \lambda_0 w_0 > z_0$.

\item[(ii)] If no such $z_0$ exists then there exists $z_0$ such that $g(z) > A_0 $ for $z < z_0$, $g(z) < A_0$ for
$z > z_0$, and the set where $\frac12 |\nabla w_0|^2 - \lambda_0 w_0 = z_0$ has positive measure. In this case
$\tilde{\theta}$ must equal $0$ where $\frac12 |\nabla w_0|^2 - \lambda_0 w_0 < z_0$ and it must
equal $1$ where $\frac12 |\nabla w_0|^2 - \lambda_0 w_0 > z_0$. (It is not fully determined by being optimal for
\eqref{optimality-of-theta-tilde}, and it could easily take values between $0$ and $1$ on the set where
$\frac12 |\nabla w_0|^2 - \lambda_0 w_0 = z_0$; this indeterminacy will not be a problem in what follows.)
\end{enumerate}
Now consider what happens to the preceding calculation when we use a different representative $w_1 = w_0 + c$. The
argument applies equally to $w_1$, except that the role of $z_0$ is played by $z_1 = z_0 - \lambda_0 c$,
since $\frac12 |\nabla w_1|^2 - \lambda_0 w_1 = \frac12 |\nabla w_0|^2 - \lambda_0 w_0 - \lambda_0 c $.

\emph{We now choose} $c = z_0/\lambda_0$, so that $z_1 = 0$ and
$$
\int_{B \setminus \overline{D}} \tilde{\theta} \left( \frac12 |\nabla w_1|^2 - \lambda_0 w_1 \right)\, dx
= \int_{B \setminus \overline{D}} \left( \frac12 |\nabla w_1|^2 - \lambda_0 w_1 \right)_+ \, dx \, ,
$$
from which it follows that
$$
L(w_1, \tilde{\theta}) \geq L(w_1, \theta) \quad \mbox{for all $\theta \in {\mathcal B}$}.
$$
We also have
$$
L(w_1, \tilde{\theta}) \leq L(w,\tilde{\theta}) \quad \mbox{for all $w \in H^1( B \setminus \overline{D})$}
$$
since $w_1$ is a representative of $\tilde{w}$. In short: $(w_1,\tilde{\theta})$ is a saddle point of $L$, viewed as a function
on $H^1(B \setminus \overline{D}) \times {\mathcal B}$. Thus, we have proved part (a) of the theorem,
with $(\overline{w},\overline{\theta}) = (w_1,\tilde{\theta})$. Part (b) is also clear from the preceding
arguments -- though it is not really necessary to check, since in general the existence of
a saddle point $(\overline{w},\overline{\theta})$ implies that the sup-inf and inf-sup are equal, and that their
common value is $L(\overline{w},\overline{\theta})$ (see e.g. Proposition 1.2 in Chapter 6 of \cite{ET}).
\end{proof}

Theorem \ref{t.round-D} showed that when $D$ is a ball and $f=\partial_{\nu_D} \psi_{d,\lambda_0}$, the
unique optimal $\Omega$ is a concentric ball. It is natural to ask whether uniqueness holds even in the larger
class of relaxed designs. The following result provides an affirmative answer -- not only when $D$ is a ball, but
also for any $D$ such that there exists an optimal (unrelaxed) ENZ shell. In addition, this
result and its proof provide a fresh perspective on the argument we used for Theorem \ref{t.round-D}.

\begin{proposition} \label{classical-solutions-of-relaxed-problem}
Let $\overline{w} \in H^1(B \setminus \overline{D})$ and $\overline{\theta} \in {\mathcal B}$ be a saddle
point for the functional $L$ defined by \eqref{definition-of-L}. Suppose furthermore that $\overline{\theta}$
solves the unrelaxed optimal design problem -- in other words that it is the
characteristic function of $\Omega_0 \setminus \overline{D}$ for some connected Lipschitz
domain $\Omega_0$ which contains $\overline{D}$ and is compactly contained in $B$. Then:
\begin{enumerate}
\item[(i)] $\Omega_0 \setminus \overline{D}$ is exactly the subset of $B \setminus \overline{D}$ where
$\frac12 |\nabla \overline{w}|^2 - \lambda_0 \overline{w} > 0$; and
\item[(ii)] for any other saddle point $(\hat{w}, \hat{\theta})$ of $L$, we have $\hat{\theta} = \overline{\theta}$
and there is a constant $c$ such that $\hat{w} = \overline{w} + c$ in $\Omega_0 \setminus \overline{D}$.
\end{enumerate}
If we assume a little more regularity -- specifically, if we assume that $\Omega_0$ is a $C^2$ domain, and that
$\nabla \overline{w} (x)$ is uniformly continuous as $x$ approaches $\partial \Omega_0$ from within $\Omega_0$ --
we can say further that
\begin{enumerate}
\item[(iii)] $\overline{w} = 0$ at $\partial \Omega_0$.
\end{enumerate}
It follows that the constant $c$ in part (ii) is actually $0$, and that $\overline{w}$ satisfies the overdetermined boundary condition
\begin{equation} \label{overdetermined-bc}
\mbox{$\partial_\nu \overline{w}= 0$ and $\overline{w} = 0$ at $\partial \Omega_0$.}
\end{equation}
\end{proposition}

\begin{remark}
Our hypothesis that $\nabla \overline{w} (x)$ be uniformly continuous as $x$ approaches $\partial \Omega_0$ from within
$\Omega_0$ follows from standard elliptic regularity results if $\partial \Omega_0$ is smooth enough.
\end{remark}

\begin{remark} \label{we-wonder}
We saw in Section \ref{subsec:round-D} that when $D$ is a ball, a concentric ball with the right area is optimal.
The proof used the associated $\overline{w}$, which vanished at the boundary of the concentric ball. (The PDE
that $\overline{w}$ solves in $\Omega_0 \setminus \overline{D}$ determines it only up to an additive constant;
however we saw in the proof of Theorem \ref{t.relaxed-design-thm} why for a saddle point of $L$ we need $\overline{w}$
to vanish -- rather than simply being constant -- on $\partial \Omega_0$.) It is not surprising that in a shape
optimization problem, the associated PDE should satisfy an overdetermined boundary condition. But we wonder
whether, for general $D$, there is really a domain $\Omega_0$ containing $\overline{D}$ for which there exists
a solution of $-\Delta w = \lambda_0$ in $\Omega_0 \setminus \overline{D}$ with $\partial_{\nu_D} w = f$ at
$\partial D$ and the overdetermined condition \eqref{overdetermined-bc} at $\partial \Omega_0$. If not, then
our optimal design problem would have no (sufficiently regular) classical solution, though it always has
a relaxed solution.
\end{remark}

\begin{proof}[Proof of Proposition \ref{classical-solutions-of-relaxed-problem}]
We have assumed that $\Omega_0$ is a connected Lipschitz domain, but we have not assumed that it
is simply connected. Thus $\Omega_0 \setminus \overline{D}$ is a bounded and connected domain in
$\R^2$ which could have finitely many ``holes.''

To get started, we collect some easy observations that follow from $(\overline{w},\overline{\theta})$
being a saddle point. The first is that $\overline{\theta}$ must actually be in $\tilde{\mathcal B}$ (defined
by \eqref{B-tilde-defn}), since otherwise $\min_{w \in H^1(B \setminus \overline{D})} L(w, \overline{\theta})$ would
be $-\infty$. Our second observation is that $-\Delta \overline{w} = \lambda_0$ in
$\Omega_0 \setminus \overline{D}$ with $\partial_\nu \overline{w} = 0$ at $\partial \Omega_0$ and
$\partial_{\nu_D} \overline{w} = f$ at $\partial D$,
since $\overline{w}$ minimizes $L(w, \overline{\theta})$ and $\overline{\theta}$ is the characteristic
function of $\Omega_0 \setminus \overline{D}$. Our third observation is that
$\frac12 |\nabla \overline{w}|^2 - \lambda_0 \overline{w}$ cannot be constant on a set
of positive measure in $\Omega_0 \setminus \overline{D}$, since
$$
\Delta \left( \frac12 |\nabla \overline{w}|^2 - \lambda_0 \overline{w} \right) = |\nabla\nabla w|^2 + \lambda_0^2 > 0 \, .
$$
(We note that, by elliptic regularity, that $\overline{w}$ is smooth in the interior of $\Omega_0 \setminus \overline{D}$.)
Our fourth observation is that
$$
\frac12 |\nabla \overline{w}|^2 - \lambda_0 \overline{w} \geq 0 \quad \mbox{in $\Omega_0 \setminus \overline{D}$, and} \quad
\frac12 |\nabla \overline{w}|^2 - \lambda_0 \overline{w} \leq 0 \quad \mbox{outside of $\Omega_0$,}
$$
since $\overline{\theta}$ maximizes
$\int_{B \setminus \overline{D}} \theta \left(\frac12 |\nabla \overline{w} |^2  - \lambda_0 \overline{w} \right) \, dx $
over all $\theta$ such that $0 \leq \theta(x) \leq 1$.

Assertion (i) of the Proposition is now easy: combining our third and fourth observations, $\Omega_0$ must agree a.e. with the set where
$\frac12 |\nabla \overline{w}|^2 - \lambda_0 \overline{w} > 0$.

For assertion (ii), we argue as we did for Theorem \ref{t.round-D}:
\begin{align}
L(\hat{w},\hat{\theta}) & \leq L(\overline{w},\hat{\theta}) \label{first-inequality} \\
 & = \int_{B \setminus \overline{D}} \hat{\theta} \left(\frac12 |\nabla \overline{w} |^2  - \lambda_0 \overline{w} \right) \, dx +
    \int_{\partial D} f \overline{w} \,  d\sh^1 \nonumber \\
 & \leq \int_{B \setminus \overline{D}} \overline{\theta} \left(\frac12 |\nabla \overline{w} |^2  - \lambda_0 \overline{w} \right) \, dx +
      \int_{\partial D} f \overline{w} \,  d\sh^1 \nonumber \\
 &  = L(\overline{w},\overline{\theta}) \, . \nonumber
\end{align}
But the saddle value is unique (see e.g. Proposition 1.2 in Chapter 6 of \cite{ET}). Therefore both inequalities in the
preceding argument are actually equalities; in particular, $\hat{\theta}$ maximizes
$\int_{B \setminus \overline{D}} \theta \left(\frac12 |\nabla \overline{w} |^2  - \lambda_0 \overline{w} \right) \, dx $
over all $\theta$ such that $0 \leq \theta(x) \leq 1$. Arguing as for part (i), it follows that $\hat{\theta} = 1$ where
$\frac12 |\nabla \overline{w}|^2 - \lambda_0 \overline{w} > 0$. But by the argument of our first observation, $\hat{\theta}$ has
the same integral as $\overline{\theta}$. So $\hat{\theta}$ must vanish outside of $\Omega_0$.

It remains to explain $\overline{w} - \hat{w}$ must be constant in $\Omega_0$. For this, we combine \eqref{first-inequality}
with the fact that $\hat{\theta} = \overline{\theta}$. Substituting $\hat{w} = \overline{w} + (\hat{w}-\overline{w})$ into the fact that
$L(\hat{w},\overline{\theta}) = L(\overline{w}, \overline{\theta})$, expanding the square, and using the stationarity of
the functional at $\overline{w}$, we conclude that
$\int_{B \setminus \overline{D}} \overline{\theta} |\nabla \hat{w} - \nabla \overline{w} |^2 \, dx = 0$. Writing this as
$\int_{\Omega_0 \setminus \overline{D}}  |\nabla \hat{w} - \nabla \overline{w} |^2 \, dx = 0$ and using that
$\Omega_0 \setminus \overline{D}$ is connected, we conclude that $\hat{w} - \overline{w}$ must be constant on this domain.

Turning now to part (iii): consider any component $\mathcal C$ of $\partial \Omega_0$ (which is now assumed to be a finite collection of $C^2$ curves). Our arguments will be local, in a vicinity of the curve $\mathcal C$. The points to one side belong to $\Omega_0 \setminus \overline{D}$, where $\overline{w}$ satisfies $\frac12 |\nabla \overline{w}|^2 - \lambda_0 \overline{w} > 0$.
Since we have assumed that $\nabla \overline{w} (x)$ is uniformly continuous as $x$ approaches the boundary from
$\Omega_0 \setminus \overline{D} $,
we can pass to the limit in the inequality and use that $\partial_\nu \overline{w} = 0$ at $\partial \Omega_0$ to conclude that
\begin{equation} \label{ineq-from-within-omega0}
\frac12 |\partial_s \overline{w}|^2 - \lambda_0 \overline{w} \geq 0 \quad \mbox{on the chosen component $\mathcal C$,}
\end{equation}
where $\partial_s$ represents the derivative tangent to the boundary.

On the other side of $\mathcal C$ we have no PDE for $\overline{w}$, however we know that
\begin{equation} \label{bulk-ineq-outside-omega0}
\frac12 |\nabla \overline{w}|^2 - \lambda_0 \overline{w} \leq 0 \quad \mbox{outside of $\Omega_0$}\, .
\end{equation}
We shall use this to show that
\begin{equation} \label{ineq-from-outside-omega0}
\frac12 |\partial_s \overline{w}|^2 - \lambda_0 \overline{w} \leq 0 \quad \mbox{on the chosen component $\mathcal C$} \, .
\end{equation}
As a first step, we now show that $\overline{w}$ is uniformly Lipschitz continous in the complement of $\Omega_0$.
Clearly $\overline{w} \geq 0$ outside of $\Omega_0$, as a consequence of \eqref{bulk-ineq-outside-omega0}. Since
$\overline{w} \in H^1(B \setminus \overline{D})$, its $L^2$ norm is bounded, so using
\eqref{bulk-ineq-outside-omega0} once again we see that
$\int_{B \setminus \overline{\Omega_0}}
 |\nabla \overline{w}|^4 \, dx < \infty$.
Since we are in two space dimensions, it follows that $\overline{w}$ uniformly bounded on
$B \setminus \overline{\Omega_0}$. Appealing to \eqref{bulk-ineq-outside-omega0} once again,
we conclude that
$$
|\nabla \overline{w}| \leq M \quad \mbox{outside $\Omega_0$,}
$$
where $M$ is an upper bound for $\sqrt{2 \lambda_0 \overline{w}}$.

Now we combine the bulk inequality \eqref{bulk-ineq-outside-omega0} with the Lipschitz estimate to get
\eqref{ineq-from-outside-omega0}. Since $\mathcal C$ is a $C^2$ curve, we may use tubular coordinates for points in
the complement of $\Omega_0$ that lie sufficiently close to $\mathcal C$. We shall use $t$ for the distance to $\mathcal C$,
and $s$ for the arclength parameter along a curve at constant distance $t$. To bound $\partial_s \overline{w}$ at some point
$(s_0,0)$ on $\mathcal C$, we shall estimate the difference quotient $\overline{w}(s_0 + \delta,0)- \overline{w}(s_0,0)$
then pass to the limit $\delta \rightarrow 0$. Clearly
$$
\overline{w}(s_0 + \delta,0)- \overline{w}(s_0,0) = [\overline{w}(s_0 + \delta,0)- \overline{w}(s_0+ \delta,t)]
 + [\overline{w}(s_0 + \delta,t)- \overline{w}(s_0,t)] + [\overline{w}(s_0,t)- \overline{w}(s_0,0)] .
$$
The first and last terms have magnitude at most $M t$. The middle term can be estimated using \eqref{bulk-ineq-outside-omega0}
(which for a.e.~$t$ holds almost everywhere in $s$ -- we naturally restrict our attention to values of $t$ with this property).
In fact,
\begin{align*}
|\overline{w}(s_0 + \delta,t)- \overline{w}(s_0,t)| & \leq \int_{s_0}^{s_0+\delta} \partial_s \overline{w}(\sigma, t) \, d\sigma \\
& \leq \delta^{1/2} \left( \int_{s_0}^{s_0+\delta} (\partial_s \overline{w})^2(\sigma, t) \, d\sigma  \right)^{1/2}
\end{align*}
From \eqref{bulk-ineq-outside-omega0} we have
$$
(\partial_s \overline{w})^2 (\sigma, t) \leq 2 \lambda_0 \overline{w}(\sigma,t) \leq 2 \lambda_0 \overline{w}(s_0,0) + O(|\delta| + |t|) \, ,
$$
so
$$
\int_{s_0}^{s_0+\delta} (\partial_s \overline{w})^2(\sigma, t) \, d\sigma \leq
\delta [ 2\lambda_0 \overline{w}(s_0, 0) + O(|\delta| + |t|) ] \, .
$$
Combining these estimates then sending $ t \rightarrow 0$ with $\delta$ held fixed, we get
$$
|\overline{w}(s_0+\delta,0) - \overline{w}(s_0, 0)| \leq
\delta \left( 2\lambda_0 \overline{w}(s_0, 0) + O(|\delta|) \right)^{1/2} \, .
$$
Dividing by $\delta$ then passing to the limit $\delta \rightarrow 0$, we conclude that
$$
|\partial_s \overline{w} (s_0,0)| \leq \left( 2 \lambda_0 \overline{w}(s_0,0) \right)^{1/2} \, .
$$
Since $(s_0, 0)$ was an arbitrary point on the chosen component $\mathcal C$, this confirms the validity of
\eqref{ineq-from-outside-omega0}. (We note that since $\overline{w} \in H^1(B \setminus \overline{D})$, its traces
on $\partial \Omega_0$ taken from inside and outside $\Omega_0$ are the same. So while \eqref{ineq-from-within-omega0}
and \eqref{ineq-from-outside-omega0} were obtained by taking limits from opposite sides of $\mathcal C$, they
estimate the same function on $\mathcal C$.)

Combining \eqref{ineq-from-within-omega0} and \eqref{ineq-from-outside-omega0}, we have shown that
$$
\frac12 |\partial_s \overline{w}|^2 - \lambda_0 \overline{w} = 0 \quad \mbox{on the chosen component $\mathcal C$}.
$$
To see that this implies $\overline{w} = 0$, we note that since $\overline{w} \geq 0$ on $\mathcal C$, the preceding relation
can be rewritten as
\begin{equation} \label{pm-relation}
\partial_s \overline{w} = \pm \sqrt{2 \lambda_0 \overline{w}} \, .
\end{equation}
Now recall our assumption for part (iii) that $\nabla \overline{w}(x)$ is uniformly continuous as $x$ approaches
$\partial \Omega_0$ from within the domain $\Omega_0$. This implies that when viewed as a function on $\mathcal C$,
$\overline{w}$ is $C^1$. Therefore the $\pm$ sign in \eqref{pm-relation} cannot change at a point where
$\overline{w} \neq 0$. If, at some point on $\mathcal C$, the function $\overline{w}$ is strictly positive and
\eqref{pm-relation} holds with a plus sign, then $\overline{w}$ must grow as $s$ increases. Similarly, if $\overline{w}$
is strictly positive and \eqref{pm-relation} holds with a minus sign, then $\overline{w}$ must grow as $s$ decreases.
Either way we reach a contradiction, since $\mathcal C$ is a \emph{closed} curve in the plane, and $\overline{w}$ is a $C^1$
function on this curve. Since this argument applies to any component of $\partial \Omega_0$, we conclude that
$\overline{w} = 0$ on the entirety of $\partial \Omega_0$, and the proof is complete.
\end{proof}

\subsection{The physical meaning of the relaxed problem} \label{subsec:relaxation}

In the previous section, our convex relaxation of the optimal design problem was introduced by replacing
a maximization over characteristic functions $\chi(x) \in \{0,1\}$ by one over densities $\theta(x) \in [0,1]$
(see \eqref{max-min-bis} and \eqref{max-min-relaxed}). We believe that the relaxed problem and the original one are
in a certain sense equivalent. Precisely: we believe that an optimal $\theta$ for the relaxed problem
is the weak limit of a maximizing sequence for the original problem. This section
explains why we believe this, though we do not have a rigorous proof.

The basic idea is simple. If $\{\chi_k(x)\}_{k=1}^\infty$ is a maximizing sequence of characteristic functions for the
original problem \eqref{max-min-bis}, then it is easy to show the existence of a subsequence converging weakly to some
function $\theta(x)$ that takes values in $[0,1]$. If the domains $\Omega_k = \{x \, : \, \chi_k(x) = 1\}$ get increasingly
complex -- for example, if they are perforated by many small holes -- then the weak limit $\theta(x)$ represents the
\emph{asymptotic density} of material at $x$ in the limit $k \rightarrow \infty$. The asymptotic performance of $\Omega_k$
depends on more than just the density -- it is also sensitive to the microstructural geometry. To avoid discussing the
microstructure explicitly, it is natural to simply assume that for each $x$, the microstructure at $x$ is optimal given
the density $\theta(x)$. We believe that this is the effect of replacing
$\chi(x) \left(\frac12 |\nabla w|^2 - \lambda_0 w \right)$ in \eqref{max-min-bis}  by
$\theta(x) \left(\frac12 |\nabla w|^2 - \lambda_0 w \right)$ in \eqref{max-min-relaxed}.

To explain the last statement, we make recourse to the theory of homogenization. For any $\eps > 0$, let
$a^\eps_k(x) = \chi_k(x) + \eps (1-\chi_k)$, and let $w^\eps_k$ solve
\begin{equation} \label{pde-for-w_k}
\begin{aligned}
- \nabla \cdot (a^\eps_k(x) \nabla w_k) - \lambda_0 \chi_k &= 0 \quad \mbox{in $B \setminus \overline{D}$}\\
\partial_{\nu_B} w_k &= 0 \quad \mbox{at $\partial B$}\\
\partial_{\nu_D} w_k &= f \quad \mbox{at $\partial D$}.
\end{aligned}
\end{equation}
(We assume here that $\chi_k$ satisfies the consistency condition for existence of $w_k$.) As $\eps \rightarrow 0$,
$\lim_{\eps \rightarrow 0} w^\eps_k$ minimizes
$$
\int_{B \setminus \overline{D}} \chi_k(x) \left( \frac12 |\nabla w|^2 - \lambda_0 w \right) \, dx +
\int_{\partial D} f w \, d \sh^1 \, ,
$$
in other words that it achieves the minimization over $w$ in \eqref{max-min-bis} (the
proof parallels that of Lemma 3.1.5 in \cite{Allaire}).
The advantage of introducing $\eps > 0$ is that $a^\eps_k(x)$ represents a mixture of two \emph{nondegenerate} materials. Our
understanding of structural optimization is rather complete in this setting: after passing to a subsequence, the
possible \emph{homogenization limits} $a^\eps_{\rm eff}$ of $a^\eps_k$ and the possible \emph{weak limits} $\theta$ of $\chi_k$
are precisely those for which $a^\eps_{\rm eff}(x)$ lies in the \emph{G-closure} of $1$ and $\eps$ with
volume fractions $\theta(x)$ and $1-\theta(x)$ respectively
(see e.g. Theorem 3.2.1 in \cite{Allaire}). The best microstructure is the one that maximizes the energy quadratic form
$\langle a^\eps_{\rm eff} \nabla w, \nabla w \rangle$ at given volume fraction. This maximum is achieved by a layered
microstructure, using layers parallel to the level lines of $w$, which achieves the well-known \emph{arithmetic-mean bound}
\begin{equation} \label{arith-mean-bound}
\langle a^\eps_{\rm eff} \nabla w, \nabla w \rangle \leq \left( \theta + \eps (1-\theta) \right) |\nabla w|^2.
\end{equation}
In summary: for the positive-epsilon analogue of our optimal design problem, the theory of homogenization provides a relaxed problem whose
minimizers are precisely the weak limits of the minimizers of the original problem. It is obtained by replacing the term
$\chi(x) |\nabla w|^2$ by the right hand side of \eqref{arith-mean-bound} and the term $\chi(x) w$ by $\theta(x) w$. In the limit
$\eps \rightarrow 0$, this procedure gives exactly our relaxed problem.

The estimates justifying the results just summarized are not uniform as $\eps \rightarrow 0$. As a result, the preceding argument
cannot be used when $\eps = 0$. In particular, it does not constitute a proof that our relaxed design problem is equivalent to the
original one.

There is a related setting where an analogous relaxation has been justified even for $\eps = 0$.
The argument goes back to \cite{Kohn-Strang-2} and it can also be found in Section 4.2 of \cite{Allaire}.
To briefly explain the idea, let us consider the optimal design problem
\begin{equation} \label{design-independent-analogue}
\sup_{\Omega \ {\rm s.t.}\, \overline{D} \subset \Omega \subset B} \ \inf_{w \in H^1(\ENZ)}
\int_{\ENZ} \left( \frac12 |\nabla w|^2 - \gamma \right) \, dx +
\int_{\partial D} g w \, d \sh^1 \, ,
\end{equation}
where $\gamma > 0$ is a constant and $g$ is a function on $\partial D$ satisfying
$\int_{\partial D} g \, d \sh^1 = 0$ (so that the min over $w$ is bounded below). This problem is very similar to \eqref{max-min-B},
but the optimal $w_*$ is now harmonic:
\begin{equation} \label{design-independent-PDE}
\begin{aligned}
\Delta w_* &= 0 \quad \mbox{in $\ENZ$}\\
\partial_{\nu_\Omega} w_* &= 0 \quad \mbox{at $\partial \Omega$}\\
\partial_{\nu_D} w_* &= g \quad \mbox{at $\partial D$}.
\end{aligned}
\end{equation}
Since \eqref{design-independent-analogue} can be written as
$$
\sup_{\Omega \ {\rm s.t.}\, \overline{D} \subset \Omega \subset B}
- \int_{\ENZ} \left( \frac12 |\nabla w_*|^2 + \gamma \right) \, dx \, ,
$$
it seeks the domain that minimizes
$\int_\Omega \left( \frac12 |\nabla w_*|^2 + \gamma \right) \, dx$. Section 4B of \cite{Kohn-Strang-2} explains how this
problem can be approached variationally. Briefly: extending $\sigma = \nabla w_*$ by zero to the entire set $B$,
letting $\sigma$ determine $\Omega = \{ x \, : \, \sigma(x) \neq 0 \}$, and using the principle of minimum
complementary energy as an
alternative representation of $\int_\Omega \frac12 |\nabla w_*|^2$, one finds that
\eqref{design-independent-analogue} is equivalent to the minimization
$$
\inf_{
\substack{{\rm div} \, \sigma = 0 \, {\rm on} \, B \setminus \overline{D}\\
\sigma \cdot \nu_D = g \, {\rm at} \, \partial D \\
\sigma \cdot \nu_B = 0 \, {\rm at} \, \partial B
}
} \int_{B \setminus \overline{D}} \Phi_\gamma (\sigma) \, dx
$$
with
$$
\Phi_\gamma (\sigma) = \left\{
\begin{array}{ll}
\frac12 |\sigma|^2 + \gamma & \mbox{if $\sigma \neq 0$}\\
0 & \mbox{if $\sigma = 0$} \, .
\end{array} \right.
$$
The relaxation of this problem is obtained by \emph{convexifying} $\Phi_\gamma$.
A homogenization-based argument similar to the one presented earlier suggests (see \cite{Kohn-Strang-2}) that
the relaxation should be obtained by replacing $\Phi_\gamma$ with
\begin{equation} \label{homogn-based-relaxation}
\min_{0 \leq \theta \leq 1} \frac{|\sigma|^2}{2 \theta} + \gamma \theta \, .
\end{equation}
This conclusion is correct, since (as verified in \cite{Kohn-Strang-2})
\eqref{homogn-based-relaxation} is equal to the convexification of $\Phi_\gamma$.

Can our relaxation of \eqref{max-min-B} be justified using an argument analogous
to the one just summarized for \eqref{design-independent-analogue}? Perhaps, however
it seems that such an argument would require substantial new ideas.

\begin{remark}
If, as we conjecture, the relaxation considered in Section \ref{subsec:convex-relaxation} is equivalent to
considering ``homogenized'' designs, then the saddle point provided by Theorem \ref{t.relaxed-design-thm} would
have $0<\overline{\theta}<1$ in any region where homogenization occurs. Reviewing the proof of that Theorem, we see that
$\overline{w}$ would need to have $\frac12 |\nabla \overline{w}|^2 - \lambda_0 \overline{w} = 0$ in such a region.
\end{remark}

\bibliographystyle{alpha}
	\bibliography{ref-june2024}
\end{document}